\documentclass[11pt]{article}
\usepackage{amssymb}
\usepackage{latexsym}
\usepackage{amsfonts}
\usepackage{amsmath}

\newtheorem{thrm}{Theorem}[section]
\newtheorem{lemma}[thrm]{Lemma}
\newtheorem{prop}[thrm]{Proposition}

\newtheorem{cor}[thrm]{Corollary}

\numberwithin{equation}{section}
\usepackage[dvips]{color}

\def\P{\mathbb{P} }
\def\Q{\mathbb{Q} }
\def\R{\mathbb{R} }
\def\N{\mathbb{N} }

\def\B{\mathcal{B} }
\def\M{\mathcal{M} }

\makeatletter
\begin{document}
\allowdisplaybreaks

\title{\Large \bf The Seneta-Heyde scaling for supercritical super-Brownian motion
\footnote{The research of this project is supported by the National Key R\&D Program of China (No. 2020YFA0712900).}
}
\author{ \bf  Haojie Hou \hspace{1mm}\hspace{1mm}
Yan-Xia Ren\footnote{The research of this author is supported by NSFC (Grant Nos. 11671017  and 11731009) and LMEQF.\hspace{1mm} } \hspace{1mm}\hspace{1mm} and \hspace{1mm}\hspace{1mm}
Renming Song\thanks{Research supported in part by a grant from the Simons
Foundation (\#429343, Renming Song).}
\hspace{1mm} }
\date{}
\maketitle

\begin{abstract}
We consider the additive martingale $W_t(\lambda)$ and the derivative martingale
$\partial W_t(\lambda)$ for one-dimensional supercritical super-Brownian motions
with general branching mechanism.
In the critical case $\lambda=\lambda_0$, we  prove that $\sqrt{t}W_t(\lambda_0)$
converges in probability
to a positive limit, which is a constant multiple of the almost sure limit
$\partial W_\infty(\lambda_0)$ of the derivative martingale
$\partial W_t(\lambda_0)$.
We also prove that,
on the survival event,
 $\limsup_{t\to\infty}\sqrt{t}W_t(\lambda_0)=\infty$
almost surely.
\end{abstract}

\medskip

\medskip


\bigskip
\begin{doublespace}

\section{Introduction}
Let $\{Z_n, n\geq 0\}$ be a supercritical Galton-Waston process with  $Z_0 = 1$  and mean $m = \mathbb{E} Z_1 \in (1, \infty).$ It is well known that $\{m^{-n} Z_n; n \geq 0\}$ is a non-negative martingale and thus converges almost surely to a limit $W$. The Kesten-Stigum theorem says that
$W$ is non-degenerate if and only if $\mathbb{E}\left[Z_1 \log Z_1\right] < \infty$.
Seneta \cite{Sen} and  Heyde \cite{Heyde} proved that if $\mathbb{E} \left[Z_1 \log Z_1\right] = \infty$,  then there exists a non-random sequence $\{c_n\}_{n\geq 0}$ such that $Z_n / c_n$ converges almost surely to a non-degenerate random variable as $n\to\infty$. This result is known as the Seneta-Heyde theorem and the sequence $\{c_n\}$ is therefore called a Seneta-Heyde norming.

A branching random walk is defined as follows. At generation 0, there is a particle
at the origin of the real line $\mathbb{R}$. At generation $n=1$, this particle dies
and splits into a finite number of offspring.
The law of the number of offspring and the positions of the offspring relative
to their parent are given by a point process $\mathcal{Z}$.
Each of these offspring  evolves independently as its parent.
Let $\mathcal{Z}_n$ denote the point process formed by the position of the particles in the  $n$-th generation.
Biggins and Kyprianou \cite{BK1, BK2} considered the non-negative martingale
$W_n(\theta): = m(\theta)^{-n}\int \exp(-\theta x)\mathcal{Z}_n (\textup{d}x)$,
which is referred to  as the additive martingale,
where
$m(\theta )=\mathbb{E}\int \exp(-\theta x)\mathcal{Z}_1 (\textup{d}x)$.
They proved that, if $m(0)>1$ and $m(\theta)<\infty$ for some  $\theta>0$, then
the limit of $W_n(\theta)$, denoted by $W(\theta)$, is non-degenerate if and only if
$\log m(\theta)-\theta m'(\theta)/m(\theta) > 0$ (supercritical)
and $\mathbb{E}\left[ W_1(\theta)\log_+ W_1(\theta) \right] < \infty$,
where $\log_+{x}:= \max\{ \log x, 0\}.$
They also showed that, when $\log m(\theta)-\theta m'(\theta)/m(\theta) > 0$ holds but $\mathbb{E}\left[ W_1(\theta)\log_+ W_1(\theta) \right] = \infty$, there exist a Seneta-Heyde norming
$\{c_n\}_{n\geq 0}$ and a non-degenerate random variable $\Delta$ such that $W_n(\theta)/c_n$ converges to $\Delta$ in probability as $n\to\infty.$

For the critical case of  $\log m(\theta)-\theta m'(\theta)/m(\theta) = 0$,
without loss of generality,
we assume that $m(\theta) = \theta =1$.
According to \cite{BK1, BK2},
the additive martingale $W_n: = W_n(1) = \int \exp(-x)\mathcal{Z}_n(x)$ converges to $0$ almost surely.
The study of the additive martingale $W_n$ in the critical case  relies on analyzing another fundamental martingale.
Under the assumption that $\mathbb{E}\left[\int x\exp(-x) \mathcal{Z}_1 (\textup{d}x)\right] = 0$, $D_n:=\int x\exp(-x) \mathcal{Z}_n (\textup{d}x)$ is a mean $0$ martingale which is referred to as the derivative martingale.
Convergence of the derivative martingale was studied by Biggins and Kyprianou \cite{BK3}.  In order to state their result, we
introduce the following integrability conditions:
\begin{equation}\label{moment1}
\sigma^2:=\mathbb{E}\left[\int x^2e^{-x}\mathcal{Z}_1 (\textup{d}x)\right]<\infty,
\end{equation}
\begin{equation}\label{moment2}
\mathbb{E}\left[\left(\int e^{-x}\mathcal{Z}_1 (\textup{d}x)\right)\log_+^2\left(\int e^{-x}\mathcal{Z}_1 (\textup{d}x)\right)\right]<\infty,
\end{equation}
\begin{equation}\label{moment3}
\mathbb{E}\left[\left(\int \left((x)_+e^{-x}\right)\mathcal{Z}_1 (\textup{d}x)\right)\log_+\left(\int \left((x)_+e^{-x}\right)\mathcal{Z}_1 (\textup{d}x)\right)\right]<\infty,
\end{equation}
Biggins and Kyprianou \cite{BK3} proved that under the assumptions \eqref{moment1}-\eqref{moment3}, $D_n$ converges almost surely to a non-degenerate non-negative limit $D_\infty$ as $n\to\infty$,
see also A\"idekon and Shi \cite[Theorem B]{AESZ}.
Hu and Shi \cite[Theorem 1.1]{HS} proved that
there exists a deterministic sequence $(a_n)_{n\geq 1}$ such that, conditioned on  survival,  $\frac{W_n}{a_n}$ converges  in distribution to some random variable $W$ with $W>0$ a.s.
It was  further proved in A\"{i}d\'ekon and Shi \cite{AESZ} that,
under the assumptions \eqref{moment1}-\eqref{moment3},
\begin{equation}\label{convinprobab-randomwalk}
 \lim_{n\to\infty}\sqrt{n}W_n =\sqrt{\frac{2}{\pi\sigma^2}} D_\infty \quad \mbox{ in probability}.
\end{equation}
They also proved that $\limsup_{n\to\infty}\sqrt{n}W_n = +\infty $ almost surely conditioned on survival.
Under the assumption that the associated random walk is in the domain of attraction
of an $\alpha$-stable law, $\alpha\in (1, 2)$, He, Liu and Zhang \cite{HLZ} proved  $n^{1/\alpha}W_n$ converges to $C D_\infty$ in probability, where $C>0$ is a constant.
For the subcritical case $\log m(\theta)-\theta m'(\theta)/m(\theta)  < 0$,
Hu and Shi \cite[Theorem 1.4]{HS} gave some convergence results for $\log W_n (\theta)$.

A branching Brownian motion (BBM) can be defined as follows.
Initially, there is a single particle at the origin.
It lives an  exponential amount of time with parameter 1.
Each particle moves according to a  Brownian motion with drift $1$
during its lifetime
and then splits into a random number, say $L$,  of new particles.
These new particles start the same process from their place of birth behaving independently
of the others. The system goes on indefinitely, unless there is no particle
at some time.
Assume that the BBM is supercritical, i.e.,  $\mathbb{E}L>1$, and $2\mathbb{E}\left[L-1\right]=1$. Let $Z_t$ be the point process formed by the position of the particles at time $t$.
The non-negative martingale
$W_t(\theta): = e^{-(\theta-1)^2t/2}\int \exp(-\theta x)Z_t (\textup{d}x)$
is called the additive martingale and plays an important role in the study of BBMs.
It is known that the limit $W(\theta)$ of $W_t(\theta)$ is non-degenerate if and only if
$|\theta|<1$ (supercritical case) and $\mathbb{E}\left[L\log_+ L\right] < \infty$, see \cite{Ch, Ne}.
Another key object for BBMs is the derivative martingale $D_t:=\int x\exp(-x) Z_t (\textup{d}x)$
in the critical case $\theta=1$.
Yang and Ren \cite{YR} proved that $D_t$ converges almost surely to a non-degenerate non-negative limit $D_\infty$ as $t\to\infty$ if and only if $\mathbb{E}\left[L\log^2_+ L\right]<\infty$, and if $\mathbb{E}\left[L\log^2_+ L\right]<\infty$ holds, $D_\infty>0$ almost surely on the event of survival.
Fluctuation of the derivative martingale $D_t$
around its  limit $D_\infty$ was given by Maillard and Pain \cite{MP}.
The analog of \eqref{convinprobab-randomwalk} is also valid for BBMs, see \cite[(1.7)]{MP}.

In this paper we consider supercritical  super-Brownian motions in $\R$.
Let ${\mathcal B}_b(\R)$ (respectively ${\mathcal B}^+(\R)$, respectively ${\mathcal B}^+_b(\R)$)  be the set of all bounded
(respectively non-negative, respectively bounded and non-negative) real-valued
Borel functions on $\R$. Let $\mathcal{M}(\mathbb{R})$ denote the space of finite Borel measures on $\mathbb{R}$. For any $f \in B_b^+ (\mathbb{R})$ and $\mu \in \mathcal{M}(\mathbb{R})$,
 we use $\langle f, \mu\rangle$ or $\mu(f)$ to denote the integral of $f$ with respect to $\mu$ whenever the integral is well-defined.
 For simplicity, we sometimes write $\Vert \mu \Vert := \langle 1, \mu \rangle$.

We will always assume that $B=\{(B_t)_{t\geq 0}; \Pi_x, x\in \R\}$
is a Brownian motion on $\R$.  Let \emph{the branching mechanism} $\psi$ be given by
 \begin{equation}
\psi(\lambda) := -\alpha \lambda + \beta \lambda^2 +\int_{(0,\infty)} \left( e^{-\lambda x} -1 +\lambda x \right) \nu (\textup{d}x), \quad \lambda\geq 0, \label{mechanism}
\end{equation}
where $ \beta \geq 0$,
$\alpha = -\psi'(0^+)$ and $\nu$ is a measure supported on $(0, \infty)$
with $\int_{(0,\infty)} (x \land x^2) \nu (\textup{d}x) < \infty.$
There exists an $\mathcal M(\R)$-valued Markov process
$X =\{(X_t)_{t\geq 0}; \mathbb P_\mu, \mu \in \mathcal M(\R)\}$ such that
$$
\mathbb P_\mu\left[e^{-  X_t(f)}\right]	= e^{- \mu(U_tf)},	\quad t\geq 0, f \in \mathcal B^+_b(\R),
$$
where $(t,x)\mapsto  U_tf(x)$ is the  unique locally bounded non-negative map on $\mathbb R_+\times \R$ such that
$$
 U_tf(x) + \Pi_x\left[\int_0^{t} \psi\left(U_{t-s} f(B_s)\right) \mathrm{d}s\right] = \Pi_x[f(B_t)], \quad t\geq 0, x\in \R.
$$
This process $X$ is known as a super-Brownian motion with branching mechanism $\psi$. For the existence of $X$ we
refer our readers to \cite{D, Dyn1993, Dyn2002} or \cite[Section 2.3]{LZ}.

The super-Brownian motion  with branching mechanism $\psi$ is called
supercritical, critical or subcritical according to $\psi'(0^+)<0$, $\psi'(0^+)=0$ or $\psi'(0^+)>0$. In this paper we concentrate on supercritical super-Brownian motions, i.e., we assume $\psi'(0^+)<0$.
We always assume that $\psi(\infty)=\infty$ which guarantees that
the event ${\cal E}:=\{\lim_{t\to\infty}\|X_t\|=0\}$ will occur with positive probability.
Let $\lambda^*$ be the largest root of the equation $\psi(\lambda)= 0$. For any $\mu\in\M(\R)$, $\mathbb{P}
_\mu({\cal E})=e^{-\lambda^*\|\mu\|}$.

In this paper we shall also assume that
\begin{equation}
\int^\infty \frac{1}{\sqrt{\int_{\lambda^*}^\xi \psi(u)\textup{d}u}} \textup{d}\xi < \infty. \label{condition1}
\end{equation}	
Under condition \eqref{condition1},
it holds that (see, for instance, \cite{KLMR})
$\mathcal{E} = \{\exists t>0 \mbox{ such that } \Vert X_t \Vert =0  \}$.

Denote by $\mathbf 0$ the null measure on $\R$.
	Write $\mathcal M^0(\R) := \mathcal M(\R)\setminus \{\mathbf 0\}$. Set $c_\lambda = - \psi'(0^+)/\lambda +\lambda /2$ and define
$$
W_t(\lambda) := e^{-\lambda c_\lambda t}\langle e^{-\lambda \cdot}, X_t \rangle, \ \ t\geq 0,\, \lambda \in \mathbb{R} .
$$
Then according to \cite{KLMR}, for any $\mu\in\M^0(\R)$, $W(\lambda): = \{W_t(\lambda): t\geq 0\}$  is a non-negative $\mathbb{P}_\mu$-martingale and thus has an almost sure limit $W_\infty(\lambda)$.
$W(\lambda)$ is called the \emph{additive martingale}.
By \cite[Theorem 2.4]{KLMR},   $W_\infty(\lambda)$ is also an $L^1(\P_\mu)$ limit if and only if $\vert\lambda\vert < \lambda_0$ and $\int_{[1,\infty)} r(\log r)\nu(\textup{d}r) < \infty$, where $\lambda_0 = \sqrt{-2\psi'(0^+)}.$

Another important martingale $\partial W(\lambda)$,
called the \emph{derivative martingale},
is defined as follows:
$$
\partial W_t(\lambda) := \langle (\lambda t + \cdot) e^{-\lambda (c_\lambda t +\cdot)}, X_t \rangle, \ \ t \geq 0.
$$
Under condition \eqref{condition1},
Kyprianou et al. \cite[Theorem 2.4]{KLMR} proved that when $|\lambda| \geq \lambda_0$,  $\partial W_t(\lambda)$ has a $\P_\mu$  almost surely  non-negative limit $\partial W_\infty(\lambda)$
for any $\mu\in\M^0(\R)$, and  when $|\lambda|  > \lambda_0$,
$\partial W_\infty(\lambda)=0$ $\P_\mu$ almost surely.
When $|\lambda| = \lambda_0$ (called the critical case),
$\partial W_\infty(\lambda)$
is almost surely positive on $\mathcal{E}^c$ if and only if
\begin{equation}
\int_{[1,\infty)} r(\log r)^2 \nu(\textup{d}r) < \infty. \label{condition2}
\end{equation}	
In this paper we  concentrate on the critical case  $|\lambda| = \lambda_0 $.
Due to symmetry, without loss of generality, we assume $\lambda = \lambda_0 $.
The derivative martingale $\partial W_t(\lambda_0)$ plays an important role
in the study of the extremal process of super-Brownian motions,
see \cite{RSZ}.

The additive martingale $W_t(\lambda_0)$ converges to $0$ as $t\to\infty$.
The goal of this paper is to find the  rate at which $W_t(\lambda_0)$ converges to $0$.
For simplicity, we write
$$
W_t := W_t(\lambda_0),\quad  \partial W_t: = \partial W_t (\lambda_0),
\quad \partial W_\infty:= \partial W_\infty(\lambda_0).
$$

Let $\{(X^{\lambda_0}_t)_{t\geq 0}; \mathbb P_\mu, \mu \in \mathcal M(\R)\}$
be a superprocess with the same branching mechanism $\psi$ in \eqref{mechanism} and with  a Brownian motion with drift $\lambda_0$ as spatial motion.  Then $\langle f, X_t^{\lambda_0} \rangle = \langle f({\lambda_0}t + \cdot), X_t\rangle$ for any $f\in\B^+_b(\R)$.
Note that $c_{\lambda_0} = \lambda_0$, we can rewrite $W_t$ and $\partial W_t$ as
$$
W_t = \langle e^{-\lambda_0\cdot}, X_t^{\lambda_0} \rangle,\quad \partial W_t = \langle \cdot e^{-\lambda_0 \cdot }, X_t^{\lambda_0} \rangle.
$$
Write $\mathbb{P}$ as a shorthand for $\mathbb{P}_{\delta_0}$.  Throughout this paper for a probability $P$,
 we will also use $P$ to denote expectation with respect to $P$.  The main results of this paper are the following two theorems:

\begin{thrm}\label{th1}
If \eqref{condition1} and \eqref{condition2} hold,
then
$$
\lim_{t\to \infty} \sqrt{t} W_t = \sqrt{\frac{2}{\pi}} \partial W_\infty \quad \mbox{ in probability with respect to }\mathbb{P}.
$$
\end{thrm}

The following result says that the above convergence in probability can not be
strengthened to almost sure convergence.

\begin{thrm}\label{th2}
If \eqref{condition1} and \eqref{condition2} hold,
then on $\mathcal{E}^c$,
\begin{equation}\label{e:th2}
\limsup_{t\to\infty} \sqrt{t}W_t = +\infty \quad \P\mbox{-almost surely}.
\end{equation}
\end{thrm}

\section{Preliminaries}

In this section, we will introduce some useful results that will be used
later.

Recall that $\{(B_t)_{t\geq 0}; \Pi_x, x\in \R\}$ is a Brownian motion.
For any $x\in\R$, we define $\tau_{x} = \inf\{t > 0: B_t=x\}$.
It is well known that $\{e^{\lambda_0B_t - \lambda_0^2 t/2}, t\geq 0\}$ is a positive $\Pi_0$-martingale with mean $1$.
We define a martingale change of measure by
\begin{equation}
\frac{\textup{d}\Pi_0^{\lambda_0}}{\textup{d} \Pi_0}\bigg\vert_{\sigma(
B_s: 0\leq s\leq t)} = e^{\lambda_0B_t - \lambda_0^2 t/2} .\label{def-Pi-lambda0}\end{equation}
Under $\Pi^{\lambda_0}_0$, $\{B_t, t\geq 0\}$ is a Brownian motion with
drift
 $\lambda_0$ staring from $0$.
 For any $y > 0,$ we define $\widetilde{\Pi}_{y}$ by
\begin{equation} \label{def-tilde-Pi}
\frac{\textup{d} \widetilde{\Pi}_{y}}{\textup{d} \Pi_0}\bigg\vert_{\sigma (B_s: s\leq t)} = \frac{y + B_t}{y}
1_{(t<\tau_{-y})}.
\end{equation}
Under $\widetilde{\Pi}_{y}$, $\{y + B_t : t \geq 0 \}$ is a Bessel-3 process
starting from $y$ and the density of $y + B_t$ is
\begin{equation}\label{Density}
f_t(x)=
\frac{x}{y \sqrt{2\pi t}}e^{-(x-y)^2/2t} (1- e^{-2xy/t})1_{\{x >0 \} }.
\end{equation}

\subsection{Branching Markov exit measures}
For any $r\ge 0$ and $x\in \R$,  let $\{(B_t)_{t\geq r}; \Pi^{\lambda_0}_{r,x}\}$ be a Brownian motion with drift $\lambda_0$ started at $x$ at time $r$. $\Pi^{\lambda_0}_{0,x}$ is the same as $\Pi^{\lambda_0}_{x}$. Let $S = \R \times [0,\infty)$, $\mathcal{B}(S)$ be the Borel $\sigma$-field on $S$ and
$\mathcal{M}(S)$ the space of finite Borel measures on $S$.
A measure $\mu \in \mathcal{M}(\R)$ is identified with its corresponding measure on $S$ concentrated on $\R\times \{0\}$. According to Dynkin \cite{E.B1.}, there exists a family of random measures $\{(X_Q, \mathbb{P}_\mu); Q\in \mathcal{S}, \mu\in \mathcal{M}(S)\}$ such that for any $Q\in \mathcal{S}$, $\mu\in \mathcal{M}(S)$ with $\textup{supp}\  \mu \subset Q$, and bounded non-negative Borel function $f(t,x)$ on $S$,
$$
\mathbb{P}_\mu \left[ \exp\left\{-\langle f, X_Q \rangle\right\}\right] = \exp\left\{-\langle V^Q_f , \mu \rangle \right\},
$$
where $V^Q_f(x,s)$ is the unique positive solution of the equation
$$
V_f^Q (x,s) + \Pi_{s,x} \int_s^\tau \psi \left(V_f^Q(B_r,r)\right)\mathrm{d} r
= \Pi_{s,x} f(B_\tau,\tau),
$$
with $\tau:= \inf\left\{ r: (B_r,r) \notin Q\right\}$.
By \cite[(1.20)]{Dyn1993}, we have the following mean formula:
\begin{equation}\label{mean1}
\mathbb{P}_\mu\langle f, X_Q \rangle=
\langle\Pi_{s,x}\left[e^{\alpha \tau} f(B_\tau,\tau)\right], \mu\rangle.
\end{equation}

For $y > 0,  t \geq 0$, we define $D_{-y}^t:=\{(x,s): -y < x , s < t\}$.
Then the random measure $X_{D_{-y}^t}^{\lambda_0}$ is concentrated on $\partial D_{-y}^t: =\left(  \{-y\} \times [0, t) \right) \cup \left([-y, +\infty]  \times   \{t\}\right)$, and for any $ \mu \in \mathcal{M}(\mathbb{R}\times[0,\infty))$ with $\textup{supp}\ \mu \subset [-y, +\infty)\times[0,t)$, and
$f\in C_b(D_{-y}^t)$  with  $f(x,s)=f(x,0)=:f(x)$ for all $s \geq 0$,
$$
\mathbb{P}_\mu \left[\exp\left\{-\langle f, X_{D_{-y}^t}^{\lambda_0}\rangle \right\}\right]
= \exp\left\{- \langle U_f^{-y,t} ( \cdot), \mu \rangle \right\},
$$
where $U_f^{-y,t}(x,s)$  is the unique positive solution of the integral  equation
\begin{equation}\label{int1}
U_f^{-y,t} (x,s) + \Pi_{s,x}^{\lambda_0} \int_s^{t\wedge \tau_{-y}} \psi \left(U_f^{-y,t}(B_r,r)\right)\mathrm{d} r = \Pi_{s,x}^{\lambda_0}[ f( B_{t\wedge\tau_{-y}})],\quad (x,s)\in \overline {D_{-y}^t},
\end{equation}
with $\overline {D_{-y}^t}$ being the closure of ${D_{-y}^t}$.
By \eqref{mean1} and the homogeneity of Brownian motion, for any $x\in \R$, we have
\begin{equation}\label{mean2}
\mathbb{P}_{\delta_x}\langle f, X_{D_{-y}^t}^{\lambda_0}\rangle
=\Pi_{x}^{\lambda_0}\left[e^{\alpha (t\wedge \tau_{-y})} f(B_{t\wedge \tau_{-y}}
	)\right].
\end{equation}
By the time homogeneity of Brownian motion with drift
$\lambda_0$,
\eqref{int1} can be written as
$$
U_f^{-y,t} (x,s) + \Pi_{x}^{\lambda_0} \int_0^{(t-s)\wedge \tau_{-y}} \psi \left(U_f^{-y,t}(B_r,r+s)\right)\mathrm{d} r = \Pi_{x}^{\lambda_0}[ f( B_{(t-s)\wedge\tau_{-y}})],\quad (x,s)\in \overline {D_{-y}^t}.
$$
Put $u_f^{-y}(x, t-s):= U_f^{-y,t} (x,s)$. The above integral equation can be written as
$$
u_f^{-y} (x,t-s) + \Pi_{x}^{\lambda_0} \int_0^{(t-s)\wedge \tau_{-y}} \psi \left(u_f^{-y}(B_r,t-r-s)\right)\mathrm{d} r = \Pi_{x}^{\lambda_0}[ f( B_{(t-s)\wedge\tau_{-y}})],\quad (x,s)\in \overline {D_{-y}^t},
$$
which is equivalent to
\begin{equation}\label{int2}
u_f^{-y} (x,s) + \Pi_{x}^{\lambda_0} \int_0^{s\wedge \tau_{-y}} \psi \left(u_f^{-y}(B_r,s-r)\right)\mathrm{d} r = \Pi_{x}^{\lambda_0}[ f( B_{s\wedge\tau_{-y}})],\quad (x,s)\in \overline {D_{-y}^t}.
\end{equation}
The special Markov property (see \cite[Theorem 1.3]{Dyn1993}, for example) implies that, for
all $D^r_{-z}\subset D^t_{-y}$
\begin{equation}\label{SMP}
\mathbb{P}_\mu  \left[\langle f, X_{D_{-y}^t}^{\lambda_0}\rangle \Big\vert \mathcal{F}_{D_{-z}^r}^{\lambda_0}\right] =\mathbb{P}_{X_{D_{-z}^r}^{\lambda_0}}
\langle f, X_{D_{-y}^t}^{\lambda_0}\rangle,
 \end{equation}
where $\mathcal{F}_{D_{-y}^t}^{\lambda_0} := \sigma \left(X_{D_{-x}^s}^{\lambda_0}: s\leq t , x \leq y \right).$

\subsection{$\mathbb{N}$-measure and spine decomposition for $X^{\lambda_0}$}
Without loss of generality, we assume that $X$ is the coordinate process on
$\mathbb D:=\{ w= (w_t)_{t\geq 0}: w \text{ is an
$\mathcal M(\R)$-valued c\`{a}dl\`{a}g function on $[0,\infty)$}\}.$
We assume that
$(\mathcal{F}_\infty, (\mathcal{F}_t)_{t\ge 0})$ is
the natural filtration on $\mathbb D$, completed as usual with the $\mathcal{F}_\infty$-measurable and $\mathbb P_\mu$-negligible sets
for every $\mu\in\mathcal{M}(\R)$.
Let $\mathbb W^+_0$ be the family of
$\mathcal M(\R)$-valued c\`{a}dl\`{a}g
functions on $(0, \infty)$ with $\mathbf 0 $ as a trap and with
$\lim_{t\downarrow0}w_t= \mathbf 0$. $\mathbb W^+_0$ can be regarded as a subset of $  \mathbb D$.

Under condition \eqref{condition1},
$\mathbb P_{\delta_x}(X_t(1)=0)>0$ for any $x\in \R$ and $t>0$, which implies that  there exists a unique family of
$\sigma$-finite measures $\{\mathbb N_x; x\in \R\}$ on
$\mathbb W^+_0$
such that
for any $\mu\in \mathcal M(\R)$, if ${\mathcal N}(\mathrm{d}w)$ is a Poisson random measure on
$\mathbb W^+_0$ with intensity measure
$$
\mathbb N_\mu(\mathrm{d}w):=\int_{\R} \mathbb N_x(\mathrm{d}w)\mu(\mathrm{d}x),
$$
then the process defined by
$$
\widehat X_0:=\mu, \quad \widehat
 X_t:=\int_{\mathbb W^+_0}w_t{\mathcal N}(\mathrm{d}w), \quad t>0,
$$
is a realization of the superprocess
$X=\{(X_t)_{t\geq 0}; \mathbb P_\mu, \mu \in \mathcal M(\R)\}$.
Furthermore,
$ \mathbb N_x(\langle f,w_t\rangle)=
\mathbb P_{\delta_x}\langle f, X_t\rangle$
and $\mathbb N_x\left[1- \exp\left\{-\langle f, w_t \rangle \right\}\right]= -\log \mathbb{P}_{\delta_x}\left[\exp\left\{-\langle f, X_t\rangle\right\}\right]$
 for any $f\in \mathcal B^+_b(\R)$
(see \cite[Theorems 8.22 and 8.23]{LZ}).
$\{\mathbb N_x; x\in \R\}$ can be regarded as measures on $\mathbb D$ carried by
$\mathbb W^+_0$,
and are called the $\mathbb N$-measures associated to $\{\mathbb P_{\delta_x}, x\in\R\}$.
Also see \cite{DyKu} for definition of $\{\mathbb N_x; x\in \R\}$.

Next, we recall an important spine decomposition for super-Brownian motions.
The spine decomposition is related to  a martingale change of measure.
 Fix $y > 0,$ define $ V_t^{-y}$ by
\begin{equation}\label{def-V}
V_t^{-y} : = \langle (y + \cdot)e^{-\lambda_0 \cdot}, X_{D_{-y}^t}^{\lambda_0} \rangle,\quad t\geq 0.
\end{equation}
From \cite{KLMR}, we know that $V_t^{-y}$ is a positive
$\mathbb{P}$-martingale with mean $y$. Define $\mathbb{Q}^{-y}$ by
\begin{equation}\label{eq:Q}
\frac{\textup{d}\mathbb{Q}^{-y}}{\textup{d}\mathbb{P}}\bigg\vert_{\mathcal{F}_t} := \frac{1}{y}V_t^{-y},
\quad t\ge 0.
\end{equation}

We say $\{(\xi_t)_{t\geq 0}, (X^{(\mathbf{n})})_{t\ge 0}, (X^{(\mathbf{m})})_{t\ge 0}, (X'_t)_{t\geq 0}; \widetilde{\mathbb{P}}^{-y}\}$
is a \emph{spine representation} of $\{(X_t)_{t\geq 0};  \mathbb{Q}^{-y}\}$ if the following are true:

(i)
The spine process is given by $\xi:=\{\xi_t, t\geq 0\}$
such that  $\{(\xi_t+\lambda_0t+y)_{t\geq 0}; \widetilde{\mathbb{P}}^{-y}\}$  is a Bessel-3 process starting from $y$.

(ii) $N$ is a Poisson process, with parameter $2\beta$, independent of $(\xi; \widetilde{\mathbb{P}}^{-y})$. Let $D^\mathbf{n}$ be the jump times of $N$
and $D^\mathbf{n}_t : = D^\mathbf{n} \cap [0,t]$.
Given $(\xi; \widetilde{\mathbb{P}}^{-y})$ and $N$, independently for each $s\in D^\mathbf{n}$, a process $\{X^{\mathbf{n},s}, \mathbb{N}_{\xi_s}\}$
is issued at the time-space  point $(s, \xi_s)$.
For $t\ge 0$, define $X_t^{(\mathbf{n})} = \sum_{s \in D_t^{\mathbf{n}}} X_{t-s}^{\mathbf{n},s}$. $X^{(\mathbf{n})}$
is referred to as the {\it continuous immigration}.

(iii) Given $(\xi; \widetilde{\mathbb{P}}^{-y})$, let $\{R_t: t\ge 0\}$ be a point process
such that the random counting measure $\sum_{t\ge 0}\delta_{(t, R_t)}$ is a Poisson random measure on $(0, \infty)\times (0, \infty)$ with intensity $\mathrm{d} tr\nu(\mathrm{d}r)$,  let
$D^\mathbf{m}$ be the projection onto the first coordinate of the atoms $\{(s_i, r_i)\}$ of this Poisson random measure and $D^\mathbf{m}_t : = D^\mathbf{m} \cap [0,t]$.
Given $\xi$ and $R$, independently for each
$s\in D^\mathbf{m}$ and $r = R_s$, a process
$\{X^{\mathbf{m},s}, \mathbb{P}_{r\delta_{\xi_s}}\}$ is issued
at the time-space point  $(s, \xi_s)$.
For $t\ge 0$, define $X_t^{(\mathbf{m})} = \sum_{s \in D_t^{\mathbf{m}}} X_{t-s}^{\mathbf{m},s}$. $X^{(\mathbf{m})}$
is referred to as the {\it discrete immigration}.

(iv) $(X', \widetilde{\mathbb{P}}^{-y} )$  is a copy of $(X, \mathbb{P})$ and $(X', \widetilde{\mathbb{P}}^{-y} )$ is independent of $\xi$,
$\{X^{\mathbf{n},s}, \mathbb{N}_{\xi_s}\}$ and $\{X^{\mathbf{m},s}, \mathbb{P}_{r\delta_{\xi_s}}\}$.

For $t\ge 0$, define
$\widetilde{X}_t = X_t' + X_t^{(\mathbf{n})}  + X_t^{(\mathbf{m})} $.
By \cite[Theorem 7.2]{KLMR},
$$
	\{(\widetilde{X}_t)_{t\geq 0}; \widetilde{\mathbb{P}}^{-y}\} \overset{d}{=} \{(X_t)_{t\geq 0}; \mathbb{Q}^{-y}\}.
$$
$\{(\widetilde{X}_t)_{t\geq 0}; \widetilde{\mathbb{P}}^{-y}\}$ is called a spine representation of $\{(X_t)_{t\geq 0}; \mathbb{Q}^{-y}\}$.

Now we give a \emph{spine representation} of
$\{(X^{\lambda_0}_t)_{t\geq 0};  \mathbb{Q}^{-y}\}$. Define
$$\xi^{\lambda_0}:=\{\xi^{\lambda_0}_t, t\ge 0\}:= \{\xi_t + \lambda_0t, t\geq 0\},$$
then
$\{\xi^{\lambda_0}_t+y, t\geq 0; \widetilde{\mathbb{P}}^{-y}\}$ is a Bessel-3 process starting from $y$.

We construct $\{(\xi^{\lambda_0}_t)_{t\geq 0}, (X^{(\mathbf{n}), \lambda_0})_{t\geq 0}, (X^{(\mathbf{m}), \lambda_0})_{t\geq 0}, ((X^{\lambda_0})'_t)_{t\geq 0}; \widetilde{\mathbb{P}}^{-y}\}$,
called a \emph{spine representation} of $\{(X^{\lambda_0}_t)_{t\geq 0}\}$, as follows:

(i) The spine is given by $\xi^{\lambda_0}=\{\xi_t + \lambda_0 t, t\geq 0\}$ such that $(\xi^{\lambda_0}+y, \widetilde{\mathbb{P}}^{-y})$ is a Bessel-3 process starting from $y$.

(ii) $Continuous\ immigration.$ Given $\xi^{\lambda_0}$,
the continuous immigration $X_{t}^{\mathbf{n},s,\lambda_0} $ immigrated  at time $s$
is defined such that
$\forall f \in B_b^+(\mathbb{R})$,
$$
\langle f, X_{t-s}^{\mathbf{n},s,\lambda_0}\rangle=\langle f(\cdot+\lambda_0 (t-s)+\lambda_0 s), X_{t-s}^{\mathbf{n},s}\rangle =\langle f(\cdot+\lambda_0 t), X_{t-s}^{\mathbf{n},s}\rangle .
$$
The almost surely countable set of the continuous immigration times in $[0,t]$ is also given by  $D^{\mathbf{n}}_t$ as in the spine decomposition of
$\{(X_t)_{t\geq 0};  \mathbb{Q}^{-y}\}$.
Define
$X_t^{(\mathbf{n}),\lambda_0}=\sum_{s \in D_t^{\mathbf{n}}}{X}_{t-s}^{\mathbf{n},s, \lambda_0}.$

(iii) $Discrete\ immigration.$ Given $\xi^{\lambda_0}$,
the discrete immigration $X_{t}^{\mathbf{m},s,\lambda_0}$ immigrated at time $s$
is defined such that
$\forall f \in B_b^+(\mathbb{R})$,
$$\langle f, X_{t-s}^{\mathbf{m},s,\lambda_0}\rangle=\langle f(\cdot+\lambda_0 (t-s)+\lambda_0 s), X_{t-s}^{\mathbf{m},s} \rangle= \langle f(\cdot+\lambda_0 t), X_{t-s}^{\mathbf{m},s} \rangle.
$$
The almost surely countable set of the discrete immigration times in $[0,t]$ is also given by $D^{\mathbf{m}}_t$ as in the spine decomposition of
$\{(X_t)_{t\geq 0};  \mathbb{Q}^{-y}\}$.
Define
$X_t^{(\mathbf{m}),\lambda_0}=\sum_{s \in D_t^{\mathbf{m}}}{X}_{t-s}^{\mathbf{m},s, \lambda_0}.$

(iv)
$\{(X^{\lambda_0})_t', t\geq 0\}$ is define such that for any
$f \in B_b^+(\mathbb{R}), \langle f, (X^{\lambda_0})_t'\rangle = \langle f(\cdot + \lambda_0 t), X_t'\rangle.$

For any $t\geq 0$, define
\begin{equation}\label{Decomposition}
\widetilde{X}_t^{\lambda_0} := (X^{\lambda_0})'_t+ X_t^{(\mathbf{n}), \lambda_0}  + X_t^{(\mathbf{m}), \lambda_0}.
\end{equation}
\begin{prop}
\begin{equation}\label{spine-decom-3}
\{(\widetilde{X}^{\lambda_0}_t)_{t\geq 0}; \widetilde{\mathbb{P}}^{-y}\} \overset{d}{=} \{(X^{\lambda_0}_t)_{t\geq 0}; \mathbb{Q}^{-y}\}.
\end{equation}
\end{prop}
\textbf{Proof:}
By the definition of $\widetilde{X}_t^{\lambda_0},  X_{t-s}^{\mathbf{n},s,\lambda_0}$ and $X_{t-s}^{\mathbf{m},s,\lambda_0},$
\begin{align*}
\langle f, \widetilde{X}_t^{\lambda_0}\rangle =&
\langle f(\cdot+\lambda_0 t), X_t'\rangle  + \sum_{s \in D_t^{\mathbf{n}}} \langle f(\cdot+\lambda_0 t), X_{t-s}^{\mathbf{n},s} \rangle + \sum_{s \in D_t^{\mathbf{m}}} \langle f(\cdot+\lambda_0 t), X_{t-s}^{\mathbf{m},s} \rangle\\
=&\langle f(\cdot+\lambda_0 t), \widetilde{X}_t \rangle.
\end{align*}
This says that $\{(\widetilde{X}^{\lambda_0}_t)_{t\geq 0}, \widetilde{\mathbb{P}}^{-y}\}$ is a shift of
$\{(\widetilde X_t)_{t\geq 0}, \widetilde{\mathbb{P}}^{-y}\}$ with constant speed $\lambda_0$.
Also note that
\begin{align*}
\Q^{-y}\left[\exp\left\{-\langle f, X^{\lambda_0}_t\rangle\right\}\right]=&\Q^{-y}\left[\exp\left\{-\langle f(\cdot+\lambda_0t),  X_t\rangle\right\}\right]=
\widetilde \P^{-y}\left[\exp\left\{-\langle f(\cdot+\lambda_0t),  \widetilde X_t\rangle\right\}\right].
\end{align*}
Thus we have
$$
\Q^{-y}\left[\exp\left\{-\langle f, X^{\lambda_0}_t\rangle\right\}\right]=\widetilde \P^{-y}\left[\exp\left\{-\langle f,  \widetilde{X}_t^{\lambda_0}\rangle\right\}\right],
$$
which says that $\{(\widetilde{X}^{\lambda_0}_t)_{t\geq 0}, \widetilde{\mathbb{P}}^{-y}\}$ and $ \{(X^{\lambda_0}_t)_{t\geq 0}, \mathbb{Q}^{-y}\}$ have the same marginal distribution. By the Markov property of both processes, we have \eqref{spine-decom-3}.
\hfill$\Box$

\subsection{Skeleton decomposition for $X$}\label{Skeleton}
In this subsection, we recall the skeleton decomposition,
which is also called the backbone decomposition in some papers, see Eckhoff et al. \cite{EKW} for an explanation of the terminologies.
This decomposition will be used in the proof of Theorem \ref{th2}.
		
Recall that $X=\{(X_t)_{t\geq0}; \P_\mu, \mu \in \mathcal{M}(\R)\}$ is a supercritical super-Brownian motion and $\mathcal{E}=\{\lim_{t\to\infty}\Vert X_t\Vert =0 \}.$
Under condition \eqref{condition1}, $\mathcal{E}=\{\Vert X_t\Vert =0 \mbox{ for some } t>0\}.$
For any $\mu \in \mathcal{M}(\R)$, we define $\P^{\mathcal{E}}_\mu$ by
$$
\P^{\mathcal{E}}_\mu(\cdot):= \P_\mu(\cdot\vert \mathcal{E}).
$$
Then by \cite[Lemma 2]{BKM}, $\{(X_t)_{t\geq 0}; \P^{\mathcal{E}}_\mu\}$ is a super-Brownian motion with branching mechanism
$$
\psi^*(\lambda):=\psi(\lambda +\lambda^*) = -\alpha^* \lambda + \beta \lambda^2 +\int_{(0,\infty)} \left( e^{-\lambda x} -1 +\lambda x \right)e^{-\lambda^* x} \nu (\textup{d}x),
$$
where
$$
\alpha^* = \alpha -2\beta\lambda^*-\int_{(0,\infty)} x\left(1-e^{-\lambda^* x}\right)\nu(\mathrm{d}x)=-\psi'(\lambda^*).
$$
We denote by $\{\N_x^{\mathcal{E}}: x\in \R \}$ the $\N$-measures associated to $\{\P_{\delta_x}^{\mathcal{E}}: x\in \R\}.$

Let $\mathcal{M}_a(\R)$  be the space of finite atomic measures on $\R$.
According to Berestycki et al. \cite{BKM},  there exists a probability
space, equipped with probability measures $\{\mathbf{P}_{(\mu, \eta)}, \mu\in \mathcal{M}(\R), \eta\in\mathcal{M}_a(\R)\}$, which carries the following processes:

(i)
$\{(Z_t)_{t\geq 0}, \mathbf{P}_{(\mu, \eta)}\}$, the skeleton,  is a
branching Brownian motion with initial configuration $\eta$,
branching rate $\psi'(\lambda^*)$, and offspring distribution with generating function
\begin{equation}\label{Generating-F}
F(s):= \frac{1}{\lambda^*\psi'(\lambda^*)}\psi\left(\lambda^*(1-s)\right)+s,\quad s\in(0,1).
\end{equation}
The law of this offspring, denoted by $\{p_n: n\geq 0\}$, satisfies $p_0 = p_1 = 0$ and for $n\geq 2$,
$$
p_n= \frac{1}{\lambda^* \psi'(\lambda^*)}\left\{\beta(\lambda^*)^21_{\{n=2\}}+(\lambda^*)^n\int_{(0,\infty)}\frac{x^n}{n!}e^{-\lambda^* x}\nu(\mathrm{d}x) \right\}.
$$
For the individuals in $Z$, we will use the classical Ulam-Harris notation. Let $\mathcal{T}^Z$ denote
the set labels realized in $Z$ and let $N_t^Z\subset\mathcal{T}^Z$ denote
 the set of individuals alive at time $t$, for $u\in N_t^Z$, we use $z_u(t)$ to denote
 the position of $u$ at time $t$. The  birth time and the death time of a particle $u$ are denoted by $b_u$ and $d_u$ respectively.

(ii) $\{(X^\mathcal{E}_t)_{t\geq 0}, \mathbf{P}_{(\mu, \eta)}\}$ is
a copy of $\{(X_t)_{t\ge 0};\P_\mu^{\mathcal{E}})$.

(iii) Three different types of immigration on $Z$: $I^{\N^\mathcal{E}}= \left\{I^{\N^\mathcal{E}}_t, t\geq 0\right\}, I^{\P^\mathcal{E}}= \left\{I^{\P^\mathcal{E}}_t, t\geq 0\right\}$ and $I^B=\left\{I^B_t, t\geq 0\right\}$, which  are independent of $X^\mathcal{E}$ and, conditioned on $Z$, are independent of each other.
The three processes are described as follows:

\begin{itemize}

\item Given $Z$, independently for each $u\in\mathcal{T}^Z$, let $N^u$ be a Poisson random measure on $(b_u, d_u]$ with intensity $2\beta$ and let $s^{1, u}_i, i=1, 2, \dots,$ be the atoms of $N^u$.
The continuous immigration $I^{\N^\mathcal{E}}$ is a measure-valued process on $\R$ such that
$$
I_t^{\N^\mathcal{E}}:= \sum_{u\in\mathcal{T}^Z}\sum_{i: s^{1,u}_i\le t}
X_{t-s^{1,u}_i}^{(1,u, i)},
$$
where $X^{(1,u, i)}$ is a measure-valued process with law
$\N_{z_u(s^{1,u}_i)}^\mathcal{E}.$

\item Given $Z$, independently for each $u\in\mathcal{T}^Z$, let $\{R^u_t: t\in (b_u, d_u]\}$ be a point process such that the random counting measure
$\sum_{t\in (b_u, d_u]}\delta_{(t, R^u_t)}$ is a Poisson random measure on $(b_u, d_u]\times (0, \infty)$ with intensity
$ \mathrm{d} t re^{-\lambda^* r} \nu(\mathrm{d}r)$
and let $\{(s^{2, u}_i, r_i): i\ge 1\}$ be the atoms of this Poisson random measure.
The discrete immigration $I^{\P^\mathcal{E}}$ is a measure-valued process on $\R$ such that
$$
I_t^{\P^\mathcal{E}}:= \sum_{u\in\mathcal{T}^Z}\sum_{i: s^{2, u}_i\le t} X_{t-s^{2, u}_i}^{(2,u,i)},
$$		
where $X^{(2,u, i)}$ is a measure-valued process with law
$\P_{r_iz_u(s^{2, u}_i)}^\mathcal{E}.$

\item
The branching point immigration $I^B$ is a measure-valued process on $\R$ such that
$$I_t^{B}:= \sum_{u\in\mathcal{T}^Z}1_{\{d_u \leq t\}} X_{t-d_u}^{(3,u)},$$
here, given $Z$, independently for each $u\in\mathcal{T}^Z$ with $d_u \leq t$, $X^{(3,u)}$ is an independent copy of $X$ issued at time $d_u$ with law $\P_{Y_u\delta_{z_u(d_u)}}$, where $Y_u$ is an independent random variable with distribution
$\pi_{O_u}(\mathrm{d}y), O_u$ is the number of the offspring of $u$ and $\{\pi_n(\mathrm{d} y), n\geq 2\}$ is a sequence of probability measures such that
$$
\pi_n(\mathrm{d}y):=\frac{1}{p_n \lambda^* \psi'(\lambda^*)}\left\{\beta(\lambda^*)^2\delta_0(\mathrm{d} y)1_{\{n=2\}} + (\lambda^*)^n\frac{y^n}{n!}e^{-\lambda^* y}\nu(\mathrm{d}y) \right\}.
$$
\end{itemize}
We define $\Lambda_t= \{\Lambda_t: t\geq 0\}$ on $\R$ by
$$
\Lambda_t := X^\mathcal{E}_t + I^{\N^\mathcal{E}}_t+I^{\P^\mathcal{E}}_t + I^B_t,\quad t\geq 0.
$$
For $\mu \in \mathcal{M}(\R)$, we denote the law of a
Poisson random measure with intensity $\lambda^* \mathrm{d} \mu$ by $\mathfrak{P}_\mu$, and define $\mathbf{P}_\mu$ by
$$
\mathbf{P}_\mu := \int \mathbf{P}_{(\mu,\eta)} \mathfrak{P}_\mu(\mathrm{d} \eta).
$$
According to \cite[Theorem 2]{BKM}, for any $\mu \in \mathcal{M}(\R)$, $\{(\Lambda_t)_{t\geq 0}; \mathbf{P}_\mu\}$ is equal in law to $\{(X_t)_{t\geq 0}; \P_\mu\}$.
The branching Brownian motion $\{Z_t, t\ge 0\}$ is referred to as the skeleton process, and $\{(\Lambda_t)_{t\geq 0};\mathbf{P}_\mu \}$ is called a \emph{skeleton decomposition} of $\{(X_t)_{t\geq 0};\P_\mu\}$.

\subsection{Properties of Brownian motion and Bessel-3 process}
Recall $B=\{(B_t)_{t\geq 0}; \Pi_x, x\in \R\}$ is a Brownian motion and $\tau_{-y} = \inf\{t > 0: B_t=-y \}$ for $y\in\R$.
\begin{lemma}\label{lemma9}	
For $ x\ge -y$,
$$
\Pi_x( t<\tau_{-y} ) = 2\int_0^{(y+x)/\sqrt{t} } \frac{1}{\sqrt{2\pi}} e^{-z^2/2}\textup{d}z, \quad t\ge 0.
$$
\end{lemma}
\textbf{Proof:} This can be easily obtained by the reflection principle of Brownian motion.
\hfill$\Box$

\begin{prop}\label{pro1} There exists a constant C such that
$$
\int_0^\infty  \Pi_z \left(B_s < x ,\ \min_{r \in [0,s]} B_r > 0 \right) \textup{d}s \leq C (1+x)(1+\min\{x,z\}),
		\quad x, z\ge 0.
$$
\end{prop}
\textbf{Proof:}
First note that, for any $h, t > 0$ and  $y \in \mathbb{R}$, we have
 \begin{equation}\label{ineq4_26}
\sup_{r \in \mathbb{R}}\Pi_y (r \leq B_t \leq r+h) =\sup_{r \in \mathbb{R}} \int_r^{r +h}\frac{1}{\sqrt{2\pi t}}
e^{-(u-y)^2/(2t)}
 \textup{d}u \leq \sup_{r \in \mathbb{R}} \int_r^{r+h} \frac{
    	    \mathrm{d}u	
    }{\sqrt{2\pi t}} = \frac{h}{\sqrt{2\pi t}}.
\end{equation}
Next, for any $0 \leq a < b$, $z \geq 0, t > 0$, by the Markov property, we have
\begin{align}\label{ineq4_3}
&\Pi_z\left(B_t \in [a, b], \min_{r\in[0,t]} B_r > 0 \right) \nonumber\\
&\leq \Pi_z \left(\min_{r \in [0, t/3]} B_r > 0   \right) \sup_{y >0} \Pi_y\left(B_{2t/3} \in [a,b], \min_{r \in [0, 2t/3]} B_r >0 \right).
\end{align}
It follows from Lemma \ref{lemma9} that
\begin{equation}\label{ineq4_25}
\Pi_z  \left(\min_{r \in [0, t/3]} B_r > 0   \right) \leq \sqrt{\frac{2}{\pi}} \frac{z}{\sqrt{t/3}} =  \sqrt{\frac{6}{\pi}} \frac{z}{\sqrt{t}}.
\end{equation}
The second term of right-hand of \eqref{ineq4_3} is bounded by
\begin{align}\label{ineq4_4}
& \Pi_y\left(B_{2t/3} \in [a,b], \min_{r \in [0, 2t/3]} B_r >0 \right)  \nonumber\\
&\leq \Pi_y \left(\min_{s \in[t/3, 2t/3]} (B_s - B_{2t/3}) > -b, B_0- B_{2t/3} \in [y-b, y-a] \right) \nonumber\\
&= \Pi_0 \left(\min_{s \in[0, t/3]} \tilde{B}_s > -b,  \tilde{B}_{2t/3} \in [y-b, y-a] \right)\nonumber\\
&\leq \Pi_0 \left(\min_{s \in[0, t/3]} \tilde{B}_s > -b\right) \sup_{v\in\mathbb{R}} \Pi_v(
\tilde{B}_{t/3}
 \in [y-b, y-a])\nonumber\\
&\leq  \sqrt{\frac{6}{\pi}} \frac{b}{\sqrt{t}} \frac{b-a}{\sqrt{2\pi t/3}} = \frac{3}{\pi}\frac{b(b-a)}{t},
\end{align}
where $\tilde{B}_s = B_{2t/3 -s} - B_{2t/3}$ is a Brownian motion
for $s \in [0,2t/3]$; we used the Markov property of $\tilde{B}$ at time $t/3$ in the second inequality of \eqref{ineq4_4}, and the last inequality of \eqref{ineq4_4} is due to \eqref{ineq4_25} and \eqref{ineq4_26}. Combining \eqref{ineq4_3}-\eqref{ineq4_4}, we obtain
\begin{equation}\label{ineq4_6}
\Pi_z\left(B_t \in [a, b], \min_{r\in[0,t]} B_r > 0 \right) \leq \sqrt{\frac{54}{\pi^3}} \frac{zb(b-a)}{\sqrt{t^3}}, \quad z\ge 0.
\end{equation}
If $x < z$, by the strong Markov property at $\tau_x$, we have
\begin{align}\label{ineq4_27}
&\int_0^\infty  \Pi_z \left(B_s < x ,\ \min_{r \in [0,s]} B_r > 0 \right) \textup{d}s   = \Pi_z \left[\int_0^\infty 1_{\{B_s < x ,\ \min_{r \in [0,s]} B_r > 0 \}} \textup{d}s \right] \nonumber\\
&\leq   \Pi_z \left[\int_{\tau_x }^\infty 1_{\{B_s < x ,\ \min_{r \in [\tau_x ,s]} B_r > 0 \}} \textup{d}s \right]
= \Pi_x \left[\int_0^\infty 1_{\{B_s < x ,\ \min_{r \in [0,s]} B_r > 0 \}} \textup{d}s \right]\nonumber\\
& = \int_0^\infty  \Pi_x \left(B_s < x ,\ \min_{r \in [0,s]} B_r > 0 \right) \textup{d}s .
\end{align}
Using \eqref{ineq4_6} and \eqref{ineq4_27}, we obtain that
\begin{align}\label{case1}
\int_0^\infty  \Pi_z \left(B_s < x ,\ \min_{r \in [0,s]} B_r > 0 \right) \textup{d}s  & \leq x^2 + \int_{x^2}^\infty  \Pi_x \left(B_s < x ,\ \min_{r \in [0,s]} B_r > 0 \right) \textup{d}s  \nonumber\\
&  \leq x^2 + \int_{x^2}^\infty  \sqrt{\frac{54}{\pi^3}} \frac{x^3}{\sqrt{s^3}}\textup{d}s \leq C_1(1+x)^2
\end{align}
for some constant $C_1>0$.
If $x \geq z$, by \eqref{ineq4_25} and \eqref{ineq4_6}, we also have
\begin{align}\label{case2}
		&\int_0^\infty  \Pi_z \left(B_s < x ,\ \min_{r \in [0,s]} B_r > 0 \right) \textup{d}s\nonumber\\
   &\leq  \int_0^{x^2} \Pi_z \left(\min_{r \in [0,s]} B_r > 0  \right)\textup{d}s   + \int_{x^2}^\infty  \Pi_z \left(B_s < x ,\ \min_{r \in [0,s]} B_r > 0 \right) \textup{d}s \nonumber \\
		 &\leq\int_0^{x^2} \sqrt{\frac{6}{\pi}}\frac{z}{\sqrt{s}}\textup{d}s +  \int_{x^2}^\infty \sqrt{\frac{54}{\pi^3}} \frac{zx^2}{\sqrt{s^3}}\textup{d}s \leq C_2(1+x)(1+z)
\end{align}
for some constant $C_2>0$. Combining \eqref{case1} and \eqref{case2}, we
arrive at the assertion of the proposition.
\hfill$\Box$

The following is a direct consequence of \cite[(3.1)]{Im}.
			
\begin{lemma}\label{lemm8}
Suppose that $\{(\eta_t)_{t\geq 0}; \widetilde\Pi_x, x\in \R_+\}$ is a Bessel-3 process.
If $F$ is a non-negative function on $C([0, t], \R)$, then
$$
\Pi_x \left[F\left(B_s, s\in [0,t]\right) 1_{\{\forall s \in [0,t], B_s > 0 \} }\right] = \widetilde{\Pi}_x \left[\frac{x}{\eta_t} F\left(\eta_s , s\in [0,t]\right) \right], \quad x\in \R_+.
$$
\end{lemma}

\begin{lemma}\label{lemma1}
If $\{(\eta_t)_{t\geq 0}; \widetilde\Pi_y, y\in \R_+\}$ is a Bessel-3 process,  then
$$
\widetilde{\Pi}_y\left[\eta^{-2}_t\right]		
\leq \frac{2}{t}, \quad t>0, \ y\ge 0.
$$
\end{lemma}
\textbf{Proof:}
Using the inequality $1-e^{-x} \leq x$ and the density of $\eta_t $ given by \eqref{Density}, we have
$$
\widetilde{\Pi}_y\left[\eta^{-2}_t\right]=
\int_0^\infty x^{-2} f_{\eta_t}(x) \textup{d}x
\leq \int_{-\infty}^\infty x^{-2} \cdot
\frac{2 x^2}{t \sqrt{2\pi t}} e^{-(x-y)^2/2t}
\textup{d}x = \frac{2}{t} .
$$
\hfill$\Box$

\begin{lemma}\label{lm2}
Suppose  that  $\{(\eta_t)_{t\geq 0}; \widetilde\Pi_y, y\in \R_+\}$ is a Bessel-3 process,
then for any event $A_t$ with  $ \lim_{t \to \infty} \widetilde{\Pi}_y(A_t) = 1$, we have
\begin{equation}\label{lemma2}
\lim_{t \to \infty} t\widetilde{\Pi}_y
\left[\eta^{-2}_t1_{A_t^c}\right]=0.
	\end{equation}
\end{lemma}
\textbf{Proof:}
For any $\varepsilon > 0,$ we have
\begin{align}
\widetilde{\Pi}_y
\left[\eta^{-2}_t1_{A_t^c}\right]
& \leq \widetilde{\Pi}_y
\left[\eta^{-2}_t1_{A_t^c}1_{\{\eta_t \geq \varepsilon \sqrt{t}\}}\right]+
\widetilde{\Pi}_y \left[\eta^{-2}_t1_{\{\eta_t < \varepsilon \sqrt{t}\} }\right]\nonumber\\
&\leq 	\widetilde{\Pi}_y(A_t^c) \cdot \frac{1}{\varepsilon^2 t} + \widetilde{\Pi}_y
\left[\eta^{-2}_t1_{\{\eta_t < \varepsilon \sqrt{t} \}}\right].
\label{step_7}
\end{align}
By the same estimate for the density of $\eta_t$ in Lemma \ref{lemma1},
\begin{align}
&\widetilde{\Pi}_y
\left[\eta^{-2}_t1_{\{\eta_t < \varepsilon \sqrt{t} \}}\right]
		=\int_0^{\varepsilon \sqrt{t}} x^{-2}f_{\eta_t}(x) \textup{d}x\nonumber\\
&\leq \frac{2}{t} \int_0^{\varepsilon \sqrt{t}}
\frac1{\sqrt{2\pi t}} e^{-(x-y)^2/2t}\textup{d}x
\leq  \frac{2}{t} \int_0^{\varepsilon \sqrt{t}} \frac{1}{\sqrt{2\pi t}} \textup{d}t = \frac{2\varepsilon}{\sqrt{2\pi}} \frac{1}{t} . \label{step_5}
\end{align}
Combining \eqref{step_7} and \eqref{step_5}, letting $t \to \infty,$ we get
$$\limsup_{t \to \infty} t\widetilde{\Pi}_y \left[\eta_t^21_{A_t^c} \right] \leq \frac{2\varepsilon}{\sqrt{2\pi}}.$$
Since $\varepsilon$ is arbitrary, we get \eqref{lemma2}.
\hfill$\Box$

\section{Proof of Theorem \ref{th1}}\label{proof}

\begin{prop}\label{prop2_5} For any $y>0$, we have
$$
\widetilde{\mathbb{P}}^{-y}\left[\xi^{\lambda_0}_t \in \textup{d}x\Big\vert \widetilde{X}_{D_{-y}^t}^{\lambda_0} \right] = \frac{e^{-\lambda_0 x}(x+y)\widetilde{X}_{D_{-y}^t}^{\lambda_0}(\textup{d}x)}{
		\widetilde V_t^{-y}},
$$
where
$$
\widetilde V_t^{-y} : = \langle (y + \cdot)e^{-\lambda_0 \cdot}, \widetilde X_{D_{-y}^t}^{\lambda_0} \rangle.
$$
\end{prop}
\textbf{Proof:}
The main idea comes from   \cite[Theorem 5.1]{KLMR}. Let $C_b^+(\partial D_{-y}^t)$ be the set of bounded non-negative continuous functions on $\partial D_{-y}^t$. We only need to show that for any $g \in C_b^+(\partial D_{-y}^t),$
\begin{equation}
\widetilde{\mathbb{P}}^{-y}\left[\exp\left\{-\theta \xi^{\lambda_0}_t-\langle g, \widetilde{X}_{D_{-y}^t}^{\lambda_0}\rangle\right\}  \right] = \widetilde{\mathbb{P}}^{-y}\left[\exp\left\{-\langle g, \widetilde{X}_{D_{-y}^t}^{\lambda_0}\rangle\right\} \frac{\langle e^{-(\lambda_0 + \theta)\cdot} (\cdot +y), \widetilde{X}_{D_{-y}^t}^{\lambda_0} \rangle}{\widetilde V_t^{-y}}\right]. \label{ineq3_4}
\end{equation}
By \eqref{spine-decom-3} and the definition \eqref{eq:Q} of $\mathbb{Q}^{-y}$,
the right hand side of \eqref{ineq3_4} is equal to
\begin{align*}
\frac{1}{y} \mathbb{P} \left[\exp\left\{-\langle g, {X}_{D_{-y}^t}^{\lambda_0}\rangle\right\} \cdot \langle e^{-(\lambda_0 + \theta)\cdot} (\cdot +y), X_{D_{-y}^t}^{\lambda_0} \rangle\right]
= -\frac{1}{y}\mathbb{P}\left[\frac{\partial}{\partial \gamma }\left[\exp\left\{-\langle g_\gamma , X_{D_{-y}^t}^{\lambda_0}\rangle\right\}\right]\bigg\vert_{\gamma = 0^+} \right]
\end{align*}
with $g_\gamma(x, t)= g(x,t) + \gamma e^{-(\lambda_0+\theta )x}(x+y).$
Interchanging the order of expectation and differentiation, we get
that
$$
\mbox{the right hand side of \eqref{ineq3_4}}
=-\frac{1}{y}\frac{\partial}{\partial \gamma }e^{- u_{g_\gamma}^{-y}(0,t)}\bigg\vert_{\gamma = 0^+} ,
$$
where $u^{-y}_{g_\gamma}$ satisfies
\eqref{int2} and $u^{-y}_{g_0} = u^{-y}_g$. Thus,
\begin{equation}\label{equal_6}
\mbox{the right hand side of \eqref{ineq3_4}}
=  \frac{1}{y}e^{- u_g^{-y}(0,t)} \frac{\partial}{\partial \gamma }u_{g_\gamma}^{-y}(0,t)\Big\vert_{\gamma = 0^+} .
\end{equation}
Let $m_g^{-y}(x,t) := \frac{\partial}{\partial \gamma }u_{g_\gamma}^{-y}(x,t)\vert_{\gamma = 0^+}$.
Replacing $f$ by $g_\gamma$ in \eqref{int2}, taking derivative with respect to $\gamma$, and
then letting $\gamma \to 0+$, we get that
$m_g^{-y}$ is the solution to the equation
$$
m_g^{-y} (x,t) + \Pi_{x}^{\lambda_0} \int_0^{t\wedge  \tau_{-y}	
 } \psi' \left(u_g^{-y}(B_r,t-r)\right)m_g^{-y}(B_r,t-r)\mathrm{d} r = \Pi_{x}^{\lambda_0}\left[ e^{-(\lambda_0+\theta )B_{t\wedge\tau_{-y}}}(B_{t\wedge\tau_{-y}}+y)\right].
 $$
Note that $B_{t	\wedge\tau_{-y}}+y=0$ when $t\geq \tau_{-y}$. The solution to the above integral equation is given by
\begin{equation}\label{equal_7}
m_g^{-y}(x,t) = \Pi_x^{\lambda_0} \left[e^{-(\lambda_0 + \theta ) B_t}(B_t +y)\exp\left\{ -\int_0^t \psi'\left(u^{-y}_g(B_{t-s}, s)\right)\textup{d}s \right\}, t<\tau_{-y} \right].
\end{equation}
By the definitions \eqref{def-Pi-lambda0} and \eqref{def-tilde-Pi}, we have
\begin{align*}
			m_g^{-y}(0,t) =& \Pi_0 \left[e^{-\frac{1}{2}\lambda^2_0t - \theta  B_t}(B_t +y)\exp\left\{ -\int_0^t \psi'\left(u^{-y}_g(B_{t-s}, s)\right)\textup{d}s \right\}, t<\tau_{-y} \right]\\
=& y\widetilde \Pi_y \left[e^{-\frac{1}{2}\lambda^2_0t - \theta  B_t}\exp\left\{ -\int_0^t \psi'\left(u^{-y}_g(B_{t-s}, s)\right)\textup{d}s \right\}\right].
\end{align*}
Using \eqref{equal_6} and \eqref{equal_7}, we have
\begin{align}\label{Laplace-xi-X}
\mbox{the right hand side of \eqref{ineq3_4}}
= 	e^{-u_g^{-y}(0,t)}\widetilde \Pi_y
	 \left[e^{-\lambda_0^2 t/2 - \theta B_t} \exp \left\{-\int_0^t \psi'\left(u^{-y}_g(B_{t-s}, s) \right)\textup{d}s  \right\} \right].
\end{align}
Next we deal with the left-hand of \eqref{ineq3_4}.
Applying Campbell's formula, we get
\begin{align}
		&	\widetilde{\mathbb{P}}^{-y}\left[\exp\left\{-\langle g, X_{D_{-y}^t}^{(\mathbf{n}), \lambda_0}\rangle\right\} \Big\vert \xi^{\lambda_0} \right] =	\widetilde{\mathbb{P}}^{-y}\left[ \exp \left\{ - \sum_{s \in D_t^{\mathbf{n}}}  \langle g, X_{D_{-y}^{t-s}}^{\mathbf{n},s, \lambda_0} \rangle
			\right\} \bigg\vert \xi^{\lambda_0} \right] \nonumber\\
&=\exp \left\{ -2\beta \int_0^t \int \left(1- \exp\left\{-\langle g, X_{D_{-y}^{ t-s}}^{\lambda_0}\rangle\right\}\right) \textup{d}\mathbb{N}_{\xi^{\lambda_0}_s} \textup{d}s \right\}  \nonumber\\
&= \exp\left\{ -2\beta \int_0^t -\log \mathbb{P}_{\delta_{\xi^{\lambda_0}_s}}\left[\exp\left\{-\langle g, X^{\lambda_0}_{D_{-y}^{t-s}} \rangle \right\} \right]\textup{d}s		\right\} \nonumber\\
&=\exp\left\{-2\beta\int_0^t u^{-y}_g(\xi^{\lambda_0}_s, t-s) \textup{d}s \right\} =\exp\left\{-2\beta\int_0^t u^{-y}_g(\xi^{\lambda_0}_{t-s}, s) \textup{d}s \right\} .\label{equal_4}
\end{align}
For $X^{(\mathbf{m}), \lambda_0},$
let $m_s:= \Vert X_{D_{-y}^0}^{\mathbf{m},s, \lambda_0} \Vert$ denote by the initial mass of the discrete immigration for $s \in D^{\mathbf{m}}$, then $\{m_s:s\geq 0\}$ is a Poisson point process
on $(0,\infty)^2$ with intensity $\mathrm{d} t  r\nu(\mathrm{d}r).$
 We similarly have
\begin{align}
&\widetilde{\mathbb{P}}^{-y}\left[\exp\left\{-\langle g, X_{D_{-y}^t}^{(\mathbf{m}), \lambda_0}\rangle \right\} \Big\vert \xi^{\lambda_0} \right]  = \widetilde{\mathbb{P}}^{-y}\left[ \exp \left\{ - \sum_{s \in
D_t^{\mathbf{m}}
}  m_s u^{-y}_g(\xi^{\lambda_0}_s, t-s)
			\right\} \bigg\vert \xi^{\lambda_0} \right] \nonumber\\
		&=	\exp\left\{ -\int_0^t \int_{(0,\infty)} \left(1-
			\exp\left\{-ru^{-y}_g(\xi^{\lambda_0}_{t-s}, s)\right\} \right)r\nu(\textup{d}r)\textup{d}s
			\right\}.\label{equal_5}
\end{align}
Combining \eqref{equal_4} and \eqref{equal_5}, we get
\begin{equation}\label{equal_12}
\widetilde{\mathbb{P}}^{-y}\left[\exp\left\{-\langle g, X_{D_{-y}^t}^{(\mathbf{n}), \lambda_0}+ X_{D_{-y}^t}^{(\mathbf{m}), \lambda_0}\rangle \right\} \Big\vert \xi^{\lambda_0} \right] = \exp\left\{-\int_0^t \left[ \psi'\left(u^{-y}_g(\xi^{\lambda_0}_{t-s}, s)\right) - \psi'(0) \right] \textup{d}s \right\}.
\end{equation}
Note that $(X^{\lambda_0})'$ is independent of $\xi$ and has the same law as $X^{\lambda_0}$, so by \eqref{equal_12},
\begin{align}\label{Laplace-xi-X'}
& 	\widetilde{\mathbb{P}}^{-y}\left[\exp\left\{-\theta \xi^{\lambda_0}_t-\langle g, \widetilde{X}_{D_{-y}^t}^{\lambda_0}\rangle\right\}  \right] \nonumber\\
&= 	\widetilde{\mathbb{P}}^{-y}\left[e^{-\theta \xi^{\lambda_0}_t}\widetilde{\mathbb{P}}^{-y}\left[\exp\left\{-\langle g,(X^{\lambda_0})'_{D_{-y}^t}+ X_{D_{-y}^t}^{(\mathbf{n}), \lambda_0}+ X_{D_{-y}^t}^{(\mathbf{m}), \lambda_0}\rangle \right\} \Big\vert \xi^{\lambda_0} \right]  \right] \nonumber\\
&=  \widetilde{\mathbb{P}}^{-y}\left[ \exp\left\{-\langle g, (X^{\lambda_0})'_{D_{-y}^t}\rangle\right\}  \right] \widetilde{\mathbb{P}}^{-y}\left[e^{-\theta \xi^{\lambda_0}_t}\widetilde{\mathbb{P}}^{-y} \left[\exp\left\{-\langle g, X_{D_{-y}^t}^{(\mathbf{n}), \lambda_0}+ X_{D_{-y}^t}^{(\mathbf{m}), \lambda_0}\rangle \right\} \Big\vert \xi^{\lambda_0} \right]  \right] \nonumber\\
&=e^{-u_g^{-y}(0,t)} \widetilde{\mathbb{P}}^{-y}\left[e^{-\theta \xi^{\lambda_0}_t}  \exp\left\{-\int_0^t \left[ \psi'\left(u^{-y}_g(\xi^{\lambda_0}_{t-s}, s)\right) - \psi'(0^+) \right] \textup{d}s \right\}\right].
\end{align}
Recall that $-\psi'(0^+) = \lambda_0^2 /2$, $\{y + B_t, t \geq 0; \widetilde{\Pi}_y \}$ is a Bessel-3 process starting from $y$ and  $\{\xi^{\lambda_0}_t+y, t\geq 0; \widetilde{\mathbb{P}}^{-y}\}$ is also a Bessel-3 process starting from $y$. Thus, by \eqref{Laplace-xi-X} and \eqref{Laplace-xi-X'}, \eqref{ineq3_4} holds.
\hfill$\Box$

For $t\geq 0$, define
\begin{equation}\label{W-T-Y}
W_t^{-y} : = \langle e^{-\lambda_0 \cdot} 1_{(-y ,\infty)}(\cdot), X_{D_{-y}^t}^{\lambda_0} \rangle
\end{equation}
and
\begin{equation}
\widetilde{W}_t^{-y} := ( W_t^{-y})' + \sum_{s \in D_t^{\mathbf{n}}}  W_{t-s}^{\mathbf{n},s,-y} + \sum_{s \in D_t^{\mathbf{m}}}   W_{t-s}^{\mathbf{m},s,-y}, \label{Decom1}
\end{equation}
where
\begin{align*}
( W_t^{-y})'&:= \langle e^{-\lambda_0 \cdot} 1_{(-y ,\infty)}(\cdot), (X^{\lambda_0})'_{D_{-y}^{t}}\rangle,\quad
W_{t-s}^{\mathbf{m},s,-y}:= \langle e^{-\lambda_0 \cdot} 1_{(-y ,\infty)}(\cdot), {X}_{D_{-y}^{t-s}}^{\mathbf{m},s, \lambda_0} \rangle,\\ W_{t-s}^{\mathbf{n},s,-y}&:= \big\langle e^{-\lambda_0 \cdot} 1_{(-y,\infty)}(\cdot), {X}_{D_{-y}^{t-s}}^{\mathbf{n},s, \lambda_0} \big\rangle.
\end{align*}
By the spine decomposition \eqref{Decomposition}, $(W_t^{-y}, t\geq 0; \mathbb{Q}^{-y})$ has the same  law as $(\widetilde W_t^{-y}, t\geq 0; \widetilde{\mathbb{P}}^{-y})$.
Recall  the definition \eqref{def-V}  of ${V}_t^{-y} $
and that $(V_t^{-y}, t\geq 0; \mathbb{Q}^{-y})$ has the same law as $(\widetilde V_t^{-y}, t\geq 0; \widetilde{\mathbb{P}}^{-y})$. Note also that
$$
\widetilde{V}_t^{-y}  = (V_t^{-y})' + \sum_{s \in D_t^{\mathbf{n}}}  V_{t-s}^{\mathbf{n},s,-y} + \sum_{s \in D_t^{\mathbf{m}}}  V_{t-s}^{\mathbf{m},s,-y},
$$
where
\begin{align*}
( V_t^{-y})' &:= \langle (y+\cdot)e^{-\lambda_0 \cdot} , (X^{\lambda_0})'_{D_{-y}^{t}}\rangle,\quad
	V_{t-s}^{\mathbf{m},s,-y} := \langle (y + \cdot)e^{-\lambda_0 \cdot}, X_{D_{-y}^{t-s}}^{\mathbf{m},s,\lambda_0} \rangle,\\
	V_{t-s}^{\mathbf{n},s,-y} &:= \langle (y + \cdot)e^{-\lambda_0 \cdot}, X_{D_{-y}^{t-s}}^{\mathbf{n},s,\lambda_0} \rangle.
\end{align*}
\begin{lemma}\label{lemma10}
For any $y > 0$ fixed, we have
$$
\lim_{t\to\infty}\sqrt{t}\,\widetilde{\mathbb{P}}^{-y} \left[ \frac{\widetilde{W}_t^{-y}}{\widetilde{W}_t^{-y} + \widetilde{V}_t^{-y}}\right]
=\sqrt{\frac2{\pi}}.
$$
\end{lemma}
\textbf{Proof:}
  First notice that
$$
 \widetilde{\mathbb{P}}^{-y} \left[ \frac{\widetilde{W}_t^{-y}}{\widetilde{V}_t^{-y}}\right] = {\mathbb{Q}^{-y}} \left[ \frac{W_t^{-y}}{V_t^{-y}}\right] = \frac{1}{y}\mathbb{P} [W_t^{-y} ].
$$
Using \eqref{mean2}, and note that $\lambda_0^2/2=\alpha$, we have that for any $f\in \mathcal{B}^+_b(\R)$,
$$
\mathbb{P}_{\delta_x}\left[\langle f, X_{D_{-y}^t}^{\lambda_0} \rangle \right] =  \Pi_x^{\lambda_0} \left[e^{\lambda_0^2 (t\land \tau_{-y})/2}f(B_{t\land \tau_{-y}}) \right].
$$
Using the above mean formula with $f(x)=e^{-\lambda_0 x}1_{(-y ,\infty)}(x)$, we obtain that
\begin{align*}
			&\widetilde{\mathbb{P}}^{-y} \left[ \frac{\widetilde{W}_t^{-y}}{\widetilde{V}_t^{-y}}\right]   =\frac{1}{y}\mathbb{P} [W_t^{-y} ]  = \frac{1}{y}\Pi_0^{\lambda_0} \left[e^{\lambda_0^2(t\land \tau_{-y})/2}e^{-\lambda_0 B_{t\land \tau_{-y}}} 1_{(-y,\infty)}( B_{t\land \tau_{-y}})\right] \\
			 &= \frac{1}{y}\Pi_0^{\lambda_0} \left[e^{\lambda_0^2 t/2}e^{-\lambda_0 B_t}1_{\{t<\tau_{-y} \}}\right]
			 = \frac{1}{y}\Pi_0 (t<\tau_{-y})  = \frac{2}{y}\int_0^{y/\sqrt{t}} \frac{1}{\sqrt{2\pi}} e^{-\frac{x^2}{2}}\textup{d}x.
\end{align*}
Thus
\begin{equation}\label{limit-W-V}
\lim_{t\to\infty}\sqrt{t}\,\widetilde{\mathbb{P}}^{-y} \left[ \frac{\widetilde{W}_t^{-y}}{ \widetilde{V}_t^{-y}}\right]=	 \sqrt{\frac2{\pi}}.
\end{equation}
To complete the proof of the lemma, it suffices to show that
$$
\limsup_{t \to \infty} \sqrt{t} \widetilde{\mathbb{P}}^{-y} \left[\frac{(\widetilde{W}_t^{-y} )^2}{(\widetilde{V}_t^{-y} +\widetilde{W}_t^{-y} )\widetilde{V}_t^{-y}} \right] = \limsup_{t \to \infty}\sqrt{t} \left\{ \widetilde{\mathbb{P}}^{-y} \left[ \frac{\widetilde{W}_t^{-y}}{ \widetilde{V}_t^{-y}}\right] -\widetilde{\mathbb{P}}^{-y} \left[ \frac{\widetilde{W}_t^{-y}}{\widetilde{W}_t^{-y} + \widetilde{V}_t^{-y}}\right]\right\} = 0.
$$
It follows from Proposition \ref{prop2_5}  that
\begin{equation}
\widetilde{\mathbb{P}}^{-y}\left[\frac{1}{
			\xi_t^{\lambda_0} + y} \Big\vert \widetilde{X}_{D_{-y}^t}^{\lambda_0} \right] = \frac{\widetilde{W}_t^{-y}}{ \widetilde{V}_t^{-y}}. \label{step_1}
\end{equation}
Under $\widetilde{\mathbb{P}}^{-y}$,
$\xi^{\lambda_0}+y$ is a Bessel-3 process starting from $y$. So by Lemma \ref{lemma1}, \eqref{step_1} and Jensen's inequality, we have
\begin{align}\label{step_1.5}
    		\widetilde{\mathbb{P}}^{-y} \left[\frac{(\widetilde{W}_t^{-y} )^2}{(\widetilde{V}_t^{-y} +\widetilde{W}_t^{-y} )\widetilde{V}_t^{-y}} \right] & \leq \widetilde{\mathbb{P}}^{-y}\left[\left(\frac{\widetilde{W}_t^{-y}}{ \widetilde{V}_t^{-y}} \right)^2 \right]  = \widetilde{\mathbb{P}}^{-y} \left[ \left(\widetilde{\mathbb{P}}^{-y}\left[\frac{1}{\xi_t^{\lambda_0} + y} \Big\vert \widetilde{X}_{D_{-y}^t}^{\lambda_0} \right]  \right)^2 \right] \nonumber\\
 & \leq \widetilde{\mathbb{P}}^{-y}\left[\left(\frac{1}{\xi_t^{\lambda_0} + y}\right)^2 \right] \leq \frac{2}{t}.
\end{align}
Therefore
$$
\sqrt{t}\widetilde{\mathbb{P}}^{-y} \left[\frac{(\widetilde{W}_t^{-y} )^2}{(\widetilde{V}_t^{-y} +\widetilde{W}_t^{-y} )\widetilde{V}_t^{-y}} \right]=o(1),\quad\mbox{ as } t\to\infty.
$$
This concludes the proof.
\hfill$\Box$\newline

Next we prove the following result:
\begin{prop}\label{convergence-L2}
\begin{equation}\label{ineq4_22}
\lim_{t\to\infty}\widetilde{\mathbb{P}}^{-y}\left[\left(\frac{\sqrt{t}\widetilde{W}_t^{-y}}{\widetilde{W}_t^{-y}+ \widetilde{V}_t^{-y}}-\sqrt{\frac2{\pi}}\right)^2\right]=0.
\end{equation}
\end{prop}

To prove \eqref{ineq4_22}, we first prove some lemmas.
Let $E_t$ be events with
$\lim_{t\to\infty}\widetilde{\mathbb{P}}^{-y}(E_t)=1$.
Combining \eqref{step_1}  and the estimate $\widetilde{\mathbb{P}}^{-y}\left[\left(\widetilde{W}_t^{-y}/ \widetilde{V}_t^{-y} \right)^2 \right] \leq\frac{2}{t}$ in \eqref{step_1.5},
we get
\begin{align}\label{step_3}
& \widetilde{\mathbb{P}}^{-y}\left[\left(\frac{\widetilde{W}_t^{-y}}{\widetilde{V}_t^{-y}+\widetilde{W}_t^{-y}} \right)^2 \right]
		\leq \widetilde{\mathbb{P}}^{-y}\left[\frac{\widetilde{W}_t^{-y}}{\widetilde{V}_t^{-y}+\widetilde{W}_t^{-y}} \frac{\widetilde{W}_t^{-y}}{\widetilde{V}_t^{-y}}  \right]\nonumber\\
&	    =	\widetilde{\mathbb{P}}^{-y}\left[\frac{\widetilde{W}_t^{-y}}{\widetilde{W}_t^{-y} + \widetilde{V}_t^{-y}} 	\widetilde{\mathbb{P}}^{-y}\left[\frac{1}{\xi_t^{\lambda_0} + y} \Big\vert \widetilde{X}_{D_{-y}^t}^{\lambda_0} \right]  \right]\nonumber\\
&=  \widetilde{\mathbb{P}}^{-y}\left[\frac{\widetilde{W}_t^{-y}}{\widetilde{W}_t^{-y}+ \widetilde{V}_t^{-y}}  \frac{1}{\xi_t^{\lambda_0} + y} \right]
		 \leq \widetilde{\mathbb{P}}^{-y}\left[\frac{\widetilde{W}_t^{-y}}{\widetilde{V}_t^{-y}}  \frac{1_{E_t^c}}{\xi_t^{\lambda_0} + y} \right] + \widetilde{\mathbb{P}}^{-y}\left[\frac{\widetilde{W}_t^{-y}}{\widetilde{W}_t^{-y} +\widetilde{V}_t^{-y}} \frac{1_{E_t}}{\xi_t^{\lambda_0} + y}   \right]
\nonumber		 \\
&\leq  \sqrt{\widetilde{\mathbb{P}}^{-y}\left[\left(\frac{\widetilde{W}_t^{-y}}{\widetilde{V}_t^{-y}} \right)^2 \right] \widetilde{\mathbb{P}}^{-y}\left[
\left(\frac{1_{E^c_t}}{\xi_t^{\lambda_0} + y} \right)^2\right] } +
 \widetilde{\mathbb{P}}^{-y}\left[\frac{\widetilde{W}_t^{-y}}{\widetilde{W}_t^{-y}  + \widetilde{V}_t^{-y}} \frac{1_{E_t}}{\xi_t^{\lambda_0} + y} \right]\nonumber\\
&\leq \sqrt{\frac{2}{t}}
\sqrt{\widetilde{\mathbb{P}}^{-y}\left[ \frac{1_{E_t^c}}{(\xi_t^{\lambda_0} + y)^2}\right] } +  \widetilde{\mathbb{P}}^{-y}\left[\frac{\widetilde{W}_t^{-y}}{\widetilde{W}_t^{-y}  +\widetilde{V}_t^{-y}} \frac{1_{E_t}}{\xi_t^{\lambda_0} + y} \right].
\end{align}
Note that, under $\widetilde \P^{-y}$, $\xi_t^{\lambda_0}+y$ is a Bessel-3 process starting from $y$.
Using Lemma \ref{lm2} and the assumption that $\widetilde{\mathbb{P}}^{-y}(E_t) \to 1$ as $t \to \infty$, we have
\begin{equation}\label{step_4}
\widetilde{\mathbb{P}}^{-y}\left[ \frac{1_{E_t^c}}{(\xi_t^{\lambda_0} + y)^2}\right] = o\left(\frac{1}{t}\right).
\end{equation}
By \eqref{step_3} and \eqref{step_4}, we conclude that
\begin{equation}\label{step_6}
\widetilde{\mathbb{P}}^{-y}\left[\left(\frac{\widetilde{W}_t^{-y}}{\widetilde{V}_t^{-y}+\widetilde{W}_t^{-y}} \right)^2 \right] \leq o\left(\frac{1}{t}\right)  + \widetilde{\mathbb{P}}^{-y}\left[\frac{\widetilde{W}_t^{-y}}{\widetilde{W}_t^{-y} +\widetilde{V}_t^{-y}} \frac{1_{E_t}}{\xi_t^{\lambda_0} + y} \right].
\end{equation}
Next, we need to construct $E_t$ such that the right-hand side of \eqref{step_6} is bounded by
$2/(\pi t) + o(1/t)$.
Let $[0, \infty)\ni t\mapsto k_t$ be a positive function such that $\lim_{t\to\infty}k_t/(\log t)^6=\infty$ and $\lim_{t\to\infty}k_t/\sqrt{t}=0$. For instance, we can take $k_t=(\log t)^7$ for large $t$.
For $t >0$ large, we define
\begin{align*}
	\widetilde{W}_t^{-y,[0, k_t)} & :=
	( W_t^{-y})'
	 +\sum_{s \in D^{\mathbf{n}}\cap[0,k_t)} W_{t-s}^{\mathbf{n},s,-y} +  \sum_{s \in D^{\mathbf{m}}\cap[0,k_t)}  W_{t-s}^{\mathbf{m},s,-y}, \\
	\widetilde{W}_t^{-y,[k_t, t]} & := \sum_{s \in D^{\mathbf{n}}\cap[k_t,t]} W_{t-s}^{\mathbf{n},s,-y} +  \sum_{s \in D^{\mathbf{m}}\cap[k_t,t]}  W_{t-s}^{\mathbf{m},s,-y}, \\
	\widetilde{V}_t^{-y,[0,k_t)} & :=
	 ( V_t^{-y})'
	 +  \sum_{s \in D^{\mathbf{n}}\cap[0,k_t)}  V_{t-s}^{\mathbf{n},s,-y} +   \sum_{s \in D^{\mathbf{m}}\cap[0,k_t)}  V_{t-s}^{\mathbf{m},s,-y}, \\
	\widetilde{V}_t^{-y,[k_t, t]} & :=  \sum_{s \in D^{\mathbf{n}}\cap[k_t,t]} V_{t-s}^{\mathbf{n},s,-y} +  \sum_{s \in D^{\mathbf{m}}\cap[k_t,t]}  V_{t-s}^{\mathbf{m},s,-y}.
\end{align*}
 Recall that $m_s= \Vert X_{D_{-y}^0}^{\mathbf{m},s, \lambda_0}\Vert$.  Define
\begin{align*}
&	E_{t,1}  := \{ k_t^{1/3} \leq \xi^{\lambda_0}_{k_t} \leq k_t \} \bigcap \left\{ \inf_{s\in [k_t, t]} \xi^{\lambda_0}_s \geq k_t^{1/6} \right\} ,
&   E_{t,2}  := \bigcap_{s\in D^\mathbf{m}\cap [k_t, t]} \left\{ m_s  \leq e^{\lambda_0 \xi^{\lambda_0}_s/2 } \right\} ,\\
&   E_{t,3} := \left\{\widetilde{V}_t^{-y,[k_t, t]} + \widetilde{W}_t^{-y,[k_t, t]} \leq \frac{1}{t^2} \right\},
&	E_t :=  E_{t,1} \cap E_{t,2}\cap E_{t,3}.
\end{align*}
\begin{lemma} \label{lemma3}
For any fixed $y>0$, it holds that
$$
\lim_{t\to\infty}\sup_{u\in[k_t^{1/3}, k_t]} \widetilde{\mathbb{P}}^{-y}\left[E_{t,2}^c\Big\vert \xi^{\lambda_0}_{k_t}=u \right] = 0.
$$
\end{lemma}
\textbf{Proof:}
First, by Campbell's formula, we have
\begin{align}\label{ineq4_1_1}
		& \widetilde{\mathbb{P}}^{-y}\left[E_{t,2}^c\Big\vert \xi^{\lambda_0}_{k_t}=u \right]
		= \widetilde{\mathbb{P}}^{-y}\left[ \bigcup_{s\in D^\mathbf{m}\cap [k_t, t]} \{m_s > e^{\lambda_0 \xi^{\lambda_0}_s/2}\}  \Big\vert \xi^{\lambda_0}_{k_t}=u \right] \nonumber\\
		&\leq   \widetilde{\mathbb{P}}^{-y}\left[\sum_{s\in D^\mathbf{m}\cap [k_t, t]} 1_{\{m_s > e^{\lambda_0 \xi^{\lambda_0}_s/2}\}} \Big\vert \xi^{\lambda_0}_{k_t}=u \right]
		\leq  \widetilde{\mathbb{P}}^{-y}\left[  \sum_{s\in D^\mathbf{m}\cap [k_t, \infty)} 1_{\{m_s > e^{\lambda_0 \xi^{\lambda_0}_s/2}\}}\bigg\vert \xi^{\lambda_0}_{k_t}=u \right]
		\nonumber\\
		&=\widetilde{\mathbb{P}}^{-y}\left[ \int_{k_t}^\infty \mathrm{d} s \int_0^\infty 1_{\{\xi_s^{\lambda_0} < 2\log r/\lambda_0  \}} r   \nu(\mathrm{d}r)\Big\vert \xi^{\lambda_0}_{k_t}=u \right].
\end{align}
Since under $\widetilde{\P}^{-y}, \xi_s^{\lambda_0} > -y$ for all $s\geq 0$, it holds that
\begin{align}\label{ineq4_1_2}
1_{\{\xi_s^{\lambda_0} < 2\log r/\lambda_0  \}} &= 1_{\{\xi_s^{\lambda_0} < 2\log r/\lambda_0  \}} \cdot 1_{\{-y < 2\log r/\lambda_0 \}} + 1_{\{\xi_s^{\lambda_0} < 2\log r/\lambda_0  \}} \cdot 1_{\{-y \geq  2\log r/\lambda_0  \}}\nonumber\\
& = 1_{\{\xi_s^{\lambda_0} < 2\log r/\lambda_0  \}} \cdot 1_{\{-y < 2\log r/\lambda_0 \}} + 1_{\{\xi_s^{\lambda_0} < 2\log r/\lambda_0 \leq -y\}} \nonumber\\
& = 1_{\{\xi_s^{\lambda_0} <2\log r/\lambda_0  \}} \cdot 1_{\{-y < 2\log r/\lambda_0 \}}.
\end{align}
Plugging \eqref{ineq4_1_2} into \eqref{ineq4_1_1} and noting
that $-y < 2\log r/\lambda_0  \Leftrightarrow r > e^{-\lambda_0 y/2}$, we get that
\begin{align}\label{ineq4_1_3}
\widetilde{\mathbb{P}}^{-y}\left[E_{t,2}^c\Big\vert \xi^{\lambda_0}_{k_t}=u \right] & \leq \widetilde{\mathbb{P}}^{-y}\left[ \int_{k_t}^\infty \mathrm{d} s \int_0^\infty 1_{\{\xi_s^{\lambda_0} < 2\log r/\lambda_0  \}} r   \nu(\mathrm{d}r)\Big\vert \xi^{\lambda_0}_{k_t}=u \right]\nonumber\\
& = \widetilde{\mathbb{P}}^{-y}\left[ \int_{k_t}^\infty \mathrm{d} s \int_{e^{-\lambda_0 y/2}}^\infty 1_{\{\xi_s^{\lambda_0} < 2\log r/\lambda_0  \}} r   \nu(\mathrm{d}r)\Big\vert \xi^{\lambda_0}_{k_t}=u \right] \nonumber\\
& = \int_{k_t}^\infty \mathrm{d} s \int_{e^{-\lambda_0 y/2}}^\infty r   \nu(\mathrm{d}r)\widetilde{\mathbb{P}}^{-y}\left[  \xi_s^{\lambda_0} < 2\log r/\lambda_0  \Big\vert \xi^{\lambda_0}_{k_t}=u \right].
\end{align}
By the Markov property, when $s \geq k_t$,
\begin{equation}\label{ineq4_1_4}
\widetilde{\mathbb{P}}^{-y}\left[  \xi_s^{\lambda_0} < 2\log r/\lambda_0  \Big\vert \xi^{\lambda_0}_{k_t}=u \right] = \widetilde{\mathbb{P}}^{-(y+u)}\left[  \xi_{s-k_t}^{\lambda_0}+u < 2\log r/\lambda_0  \right].
\end{equation}
So \eqref{ineq4_1_3} and \eqref{ineq4_1_4} yield that
\begin{align}\label{ineq4_1_5}
\widetilde{\mathbb{P}}^{-y}\left[E_{t,2}^c\Big\vert \xi^{\lambda_0}_{k_t}=u \right] &\leq  \int_{k_t}^\infty \mathrm{d} s \int_{e^{-\lambda_0 y/2}}^\infty r   \nu(\mathrm{d}r)\widetilde{\mathbb{P}}^{-(y+u)}\left[  \xi_{s-k_t}^{\lambda_0}+u < 2\log r/\lambda_0   \right]\nonumber\\
& = \int_{0}^\infty \mathrm{d} s \int_{e^{-\lambda_0 y/2}}^\infty r   \nu(\mathrm{d}r)\widetilde{\mathbb{P}}^{-(y+u)}\left[  \xi_{s}^{\lambda_0}+u < 2\log r/\lambda_0  \right].
\end{align}
Now by lemma \ref{lemm8} and Proposition \ref{pro1}, \eqref{ineq4_1_5} is bounded above by
\begin{align}\label{step_8}
&\widetilde{\mathbb{P}}^{-y}\left[E_{t,2}^c\Big\vert \xi^{\lambda_0}_{k_t}=u \right]  \leq \int_{0}^\infty \mathrm{d} s \int_{e^{-\lambda_0 y/2}}^\infty r   \nu(\mathrm{d}r)\widetilde{\mathbb{P}}^{-(y+u)}\left[  \xi_{s}^{\lambda_0}+u < 2\log r/\lambda_0  \right] \nonumber\\
&= \int_{e^{-\lambda_0y/2}}^\infty r \nu(\textup{d}r) \int_0^\infty \textup{d}s \frac{1}{u+y} \Pi_{u+y} \left(B_s 1_{\{B_s < y+2\log r/\lambda_0  , s<\tau_0 \}} \right) \nonumber\\
&\leq	\int_{e^{-\lambda_0y/2}}^\infty r \nu(\textup{d}r) \int_0^\infty \textup{d}s \frac{y+ 2\log r/\lambda_0  }{u+y} \Pi_{u+y} \left(B_s < y + 2\log r/\lambda_0  , s<\tau_0 \right)\nonumber\\
&\leq  \frac{C}{u+y} \int_{e^{-\lambda_0y/2}}^\infty r  (1+ y+2\log r/\lambda_0  )^2 \left(1+ \min\{y+2\log r/\lambda_0  , u+y  \} \right)\nu(\textup{d}r).
\end{align}
For any fixed $\varepsilon > 0$,
note that $2\log r/\lambda_0 \leq  \varepsilon u  \Longleftrightarrow r \leq  e^{\varepsilon\lambda_0u /2}.$
We suppose that $t$ is large enough
such that for any $u \in [k_t^{1/3}, k_t]$, $u+y >1$ and $1+\varepsilon u + y \leq 2\varepsilon (u +y).$ Thus,
\begin{align}\label{step_9}
\widetilde{\mathbb{P}}^{-y}\left[E_{t,2}^c\Big\vert \xi^{\lambda_0}_{k_t}=u \right]
\leq & \frac{C}{u+y} \int_{e^{-\lambda_0y/2}}^{e^{\varepsilon\lambda_0u /2}} r(1+ y + 2\log r/\lambda_0 )^2 (1+y+ 2\log r/\lambda_0 )\nu (\textup{d}r)
 \nonumber\\
 & + \frac{C(1+u+y)}{u+y}\int_{e^{\varepsilon\lambda_0u /2}} ^\infty  r(1+ y+2\log r/\lambda_0 )^2 \nu (\textup{d}r)
\nonumber \\
 \leq & \frac{C}{u+y} \int_{e^{-\lambda_0y/2}}^{e^{\varepsilon\lambda_0u /2}} r(1+ y+2\log r/\lambda_0 )^2  \left(1+ y +\varepsilon u  \right) \nu (\textup{d}r)
	\nonumber\\
& + \frac{C(1+u+y)}{u+y}\int_{ e^{\varepsilon\lambda_0u /2}}^{\infty}   r(1+ y+2\log r/\lambda_0 )^2  \nu (\textup{d}r)
\nonumber \\
 \leq &2 C \varepsilon \int_{e^{-\lambda_0y/2}}^{\infty} r(1+ y+2\log r/\lambda_0 )^2 \nu (\textup{d}r)
 \nonumber \\
& +2 C \int_{e^{\varepsilon\lambda_0 k_t^{1/3} /2}} ^\infty  r(1+ y+2\log r/\lambda_0 )^2 \nu (\textup{d}r) .
\end{align}
Using condition \eqref{condition1}
and taking $t \to \infty$, \eqref{step_9} yields that
$$
\limsup_{t \to \infty} \sup_{u \in[k_t^{1/3}, k_t]} \widetilde{\mathbb{P}}^{-y}\left[E_{t,2}^c\Big\vert \xi^{\lambda_0}_{k_t}=u \right]   \leq C \varepsilon \int_{e^{-\lambda_0y/2}}^{\infty} r(1+ y+2\log r/\lambda_0 )^2 \nu (\textup{d}r) .
$$
Since $\varepsilon $ is arbitrary, the desired assertion is valid.
\hfill$\Box$

\begin{lemma}\label{lemma4}
For any fixed $y > 0$, there exist  constants $T,
C'>0$ such that for any $t\geq T$,
$$
\widetilde{\mathbb{P}}^{-y}\left[E_{t,1}\cap E_{t,2} \cap E_{t,3}^c \big\vert \xi^{\lambda_0} \right]\leq \frac{
		C'}{t},\quad \widetilde{\mathbb{P}}^{-y} \mbox{-a.s.}
$$
\end{lemma}
\textbf{Proof:}
Recall that $W_t^{-y}$ is defined in \eqref{W-T-Y}.
	 Define $\mathcal{W}_t^{-y}$ by
	 $$\mathcal{W}_t^{-y} : = \langle e^{-\lambda_0 \cdot}, X_{D_{-y}^t}^{\lambda_0} \rangle.$$
By \eqref{mean2}, for any $t,r > 0$ and $z \geq -y$,  $\mathbb{P}_{r\delta_z}\left[\mathcal{W}_t^{-y} \right]= re^{-\lambda_0 z}$, which does not depend on $t$.
By this and the special Markov property \eqref{SMP}, we see that $\mathcal{W}_t^{-y}$ is a non-negative $\mathbb{P}_{r\delta_z}$-martingale.
Note that $W_t^{-y} \leq \mathcal{W}_t^{-y}$.
Similarly to \eqref{Decom1},we define
$$
\mathcal{W}_{t-s}^{\mathbf{m},s,-y}:= \langle e^{-\lambda_0 \cdot}, {X}_{D_{-y}^{t-s}}^{\mathbf{m},s, \lambda_0} \rangle,\quad \mathcal{W}_{t-s}^{\mathbf{n},s,-y}:= \langle e^{-\lambda_0 \cdot} , {X}_{D_{-y}^{t-s}}^{\mathbf{n},s, \lambda_0} \rangle.
$$
Because	$E_{t,1} \in \sigma(\xi_t: t \geq 0)$,
by the martingale property of $\mathcal{W}_t^{-y}$,
also by the definition of $D^{\mathbf{n}}$, we obtain that
\begin{align}
&\widetilde{\mathbb{P}}^{-y}\left[1_{E_{t,1}} \sum_{s \in D^{\mathbf{n}}\cap[k_t,t]}  W_{t-s}^{\mathbf{n},s,-y}   \bigg\vert \xi^{\lambda_0} \right]  \leq  \widetilde{\mathbb{P}}^{-y}\left[1_{E_{t,1}} \sum_{s \in D^{\mathbf{n}}\cap[k_t,t]}  \mathcal{W}_{t-s}^{\mathbf{n},s,-y}   \bigg\vert \xi^{\lambda_0} \right] \nonumber\\
&= 2\beta 1_{E_{t,1}} \int_{k_t}^t  \mathbb{N}_{\xi_s^{\lambda_0}}\left(\mathcal{W}_{t-s}^{\mathbf{n},s,-y}\Big\vert \xi^{\lambda_0}\right)\textup{d}s   = 2\beta 1_{E_{t,1}} \int_{k_t}^t \mathbb{P}_{\delta_{\xi_s^{\lambda_0}}} \left(\mathcal{W}_{t-s}^{\mathbf{n},s,-y}\Big\vert \xi^{\lambda_0}\right)\textup{d}s \nonumber\\
&= 2\beta 1_{E_{t,1}} \int_{k_t}^t e^{-\lambda_0 \xi^{\lambda_0}_s} \textup{d}s
\leq 2\beta t e^{-\lambda_0 k_t^{1/6}}
\leq 2\beta t e^{-\lambda_0 k_t^{1/6}/4},
\label{ineq4_5}
\end{align}
where the second to the last inequality of \eqref{ineq4_5} holds
because on $E_{t,1}$ we have $\xi_s \geq k_t^{1/6}$ for all $k_t \leq s \leq t$.
Next, for $ s\in D^{\mathbf{m}}$
and recall that $m_s  = \Vert X_{D_{-y}^0}^{\mathbf{m},s} \Vert$,
by the martingale property of $\mathcal{W}_t^{-y}$,
\begin{align}\label{ineq4_18}
&\widetilde{\mathbb{P}}^{-y}\left[1_{E_{t,1}\cap E_{t,2}} \sum_{s \in D^{\mathbf{m}}\cap[k_t,t]}  W_{t-s}^{\mathbf{m},s,-y}   \bigg\vert \xi^{\lambda_0}, \mathbf{m} \right]
  \leq 	\widetilde{\mathbb{P}}^{-y}\left[1_{E_{t,1}\cap E_{t,2}} \sum_{s \in D^{\mathbf{m}}\cap[k_t,t]}  \mathcal{W}_{t-s}^{\mathbf{m},s,-y}   \bigg\vert \xi^{\lambda_0}, \mathbf{m} \right]
\nonumber		\\
&=  1_{E_{t,1}\cap E_{t,2}} \sum_{s \in D^{\mathbf{m}}\cap[k_t,t]}  \mathbb{P}_{m_s\delta_{\xi_s^{\lambda_0}}}\left(\mathcal{W}_{t-s}^{\mathbf{m},s,-y} \Big\vert \xi^{\lambda_0},\mathbf{m}	\right)
  = 1_{E_{t,1}\cap E_{t,2}} \sum_{s \in D^{\mathbf{m}}\cap[k_t,t]} e^{-\lambda_0 \xi^{\lambda_0}_s} m_s \nonumber\\
& \leq  1_{E_{t,1}} \sum_{s \in D^{\mathbf{m}}\cap[k_t,t]}e^{-\lambda_0 \xi^{\lambda_0}_s/2} 1_{\{m_s > 1\}} + 1_{E_{t,1}} \sum_{s \in D^{\mathbf{m}}\cap[k_t,t]}e^{-\lambda_0 \xi^{\lambda_0}_s}m_s 1_{\{m_s \leq 1\}}
\nonumber\\
 &\leq   e^{-\lambda_0 k_t^{1/6}/2} \sum_{s \in D^{\mathbf{m}}\cap[k_t,t]}1_{\{m_s >1 \}} + e^{-\lambda_0 k_t^{1/6}} \sum_{s \in D^{\mathbf{m}}\cap[k_t,t]}m_s1_{\{m_s \leq 1 \}}.
\end{align}
Taking expectation with respect to $\mathbf{m}$ in \eqref{ineq4_18}, we get
\begin{align}\label{ineq4_19}
&\widetilde{\mathbb{P}}^{-y}\left[1_{E_{t,1}\cap E_{t,2}} \sum_{s \in D^{\mathbf{m}}\cap[k_t,t]}  W_{t-s}^{\mathbf{m},s,-y}   \bigg\vert \xi^{\lambda_0} \right] \nonumber\\
&\leq e^{-\lambda_0 k_t^{1/6}/2} \widetilde{\mathbb{P}}^{-y} \left[ \sum_{s \in D^{\mathbf{m}}\cap[k_t,t]}1_{\{m_s >1 \}} \bigg\vert \xi^{\lambda_0} \right]  + e^{-\lambda_0 k_t^{1/6}} \widetilde{\mathbb{P}}^{-y} \left[ \sum_{s \in D^{\mathbf{m}}\cap[k_t,t]}m_s1_{\{m_s \leq 1 \}} \bigg\vert \xi^{\lambda_0} \right]\nonumber\\
&= e^{-\lambda_0 k_t^{1/6}/2} \int_{k_t}^t \textup{d}s\int_1^\infty r \nu(\textup{d}r) + e^{-\lambda_0 k_t^{1/6}} \int_{k_t}^t \textup{d}s\int_0^1 r^2 \nu(\textup{d}r)\nonumber\\
&\leq  t  e^{-\lambda_0 k_t^{1/6}/2}\int_1^\infty r \nu(\textup{d}r)+ t e^{-\lambda_0 k_t^{1/6}} \int_0^1 r^2 \nu(\textup{d}r)
\leq C_3 t  e^{-\lambda_0 k_t^{1/6}/4}
\end{align}
for some constant $C_3$.
Similarly, for large $t$ such that
for all $u\geq k_t^{1/3}$, $(y+u)\leq e^{\lambda_0 u/4}$, we have
\begin{align}\label{ineq4_20}
&\widetilde{\mathbb{P}}^{-y} \left[ 1_{E_{t,1}}\sum_{s \in D^{\mathbf{n}}\cap[k_t,t]}  V_{t-s}^{\mathbf{n},s,-y} \bigg\vert \xi^{\lambda_0} \right]  = 2\beta 1_{E_{t,1}} \int_{k_t}^t  \mathbb{N}_{\xi_s^{\lambda_0}}\left(V_{t-s}^{\mathbf{n},s,-y}\Big\vert \xi^{\lambda_0}\right)\textup{d}s \nonumber\\
&= 2\beta 1_{E_{t,1}} \int_{k_t}^t  \mathbb{P}_{\delta_{\xi_s^{\lambda_0}}} \left(V_{t-s}^{\mathbf{n},s,-y}\Big\vert \xi^{\lambda_0}\right)\textup{d}s = 2\beta 1_{E_{t,1}} \int_{k_t}^t e^{-\lambda_0 \xi^{\lambda_0}_s}
				 (y+ \xi_s^{\lambda_0}) \textup{d}s\nonumber\\
& \leq  2\beta t e^{-3\lambda_0 k_t^{1/6}/4} \leq 2\beta t e^{-\lambda_0 k_t^{1/6}/4}.
\end{align}
For large $t$ such that
for all $u\geq k_t^{1/3}$, $(y+u)\leq e^{\lambda_0 u/4}$, we also have
\begin{align}\label{ineq4_7}
&\widetilde{\mathbb{P}}^{-y}\left[1_{E_{t,1}\cap E_{t,2}} \sum_{s \in D^{\mathbf{m}}\cap[k_t,t]}  V_{t-s}^{\mathbf{m},s,-y}  \bigg\vert \xi^{\lambda_0} ,	\mathbf{m} \right]
   =  1_{E_{t,1}\cap E_{t,2}} \sum_{s \in D^{\mathbf{m}}\cap[k_t,t]}  \mathbb{P}_{m_s\delta_{\xi_s^{\lambda_0}}}\left(V_{t-s}^{\mathbf{m},s,-y} \Big\vert \xi^{\lambda_0}   ,\mathbf{m}    \right)\nonumber\\
 &= 1_{E_{t,1}\cap E_{t,2}} \sum_{s \in D^{\mathbf{m}}\cap[k_t,t]} e^{-\lambda_0 \xi^{\lambda_0}_s}
   (y + \xi_s^{\lambda_0})m_s
      \leq 1_{E_{t,1}\cap E_{t,2}} \sum_{s \in D^{\mathbf{m}}\cap[k_t,t]} e^{-3\lambda_0 \xi^{\lambda_0}_s/4} m_s \nonumber\\
&\leq  e^{-\lambda_0 k_t^{1/6}/4} \sum_{s \in D^{\mathbf{m}}\cap[k_t,t]}1_{\{m_s >1 \}} + e^{-3\lambda_0 k_t^{1/6}/4} \sum_{s \in D^{\mathbf{m}}\cap[k_t,t]}m_s1_{\{m_s \leq 1 \}}.
\end{align}
Taking expectation with respect to $\mathbf{m}$
 in \eqref{ineq4_7}, we obtain that for some constant $C_4$,
\begin{align}\label{ineq4_21}
	\widetilde{\mathbb{P}}^{-y}\left[1_{E_{t,1}\cap E_{t,2}} \sum_{s \in D^{\mathbf{m}}\cap[k_t,t]} V_{t-s}^{\mathbf{m},s,-y}   \bigg\vert \xi^{\lambda_0}  \right]
	&	\leq t  e^{-\lambda_0 k_t^{1/6}/4}\int_1^\infty r \nu(\textup{d}r)+ t e^{-3\lambda_0 k_t^{1/6}/4} \int_0^1 r^2 \nu(\textup{d}r) \nonumber\\
	&\leq C_4 t  e^{-\lambda_0 k_t^{1/6}/4}.
\end{align}
Combining \eqref{ineq4_5}, \eqref{ineq4_19}, \eqref{ineq4_20} and \eqref{ineq4_21}, we get that
$$
\widetilde{\mathbb{P}}^{-y}\left[1_{E_{t,1}\cap E_{t,2}} \left(\widetilde{V}_t^{-y,[k_t, t]} + \widetilde{W}_t^{-y,[k_t, t]} \right)\bigg\vert \xi^{\lambda_0}  \right]
		\leq (C_3 + C_4 + 4\beta) t e^{-\lambda_0 k_t^{1/6}/4}.
$$
On $E_{t,3}^c $ we have $\widetilde{V}_t^{-y,[k_t, t]} + \widetilde{W}_t^{-y,[k_t, t]} > 1/t^2$.
Then for $t$ large enough such that $k_t^{1/6} > 16\log t / \lambda_0,$ we have
\begin{align*}
		&\widetilde{\mathbb{P}}^{-y} \left[1_{E_{t,1 \cap E_{t,2}\cap E_{t,3}^c} } \big\vert \xi^{\lambda_0} \right]   \leq t^2	\widetilde{\mathbb{P}}^{-y}\left[1_{E_{t,1}\cap E_{t,2}} \left(\widetilde{V}_t^{-y,[k_t, t]} + \widetilde{W}_t^{-y,[k_t, t]} \right)\bigg\vert \xi^{\lambda_0}  \right]  \\
		 &\leq (C_3 + C_4 + 4\beta) t^3 e^{-\lambda_0 k_t^{1/6}/4} \leq (C_3 + C_4 + 4\beta) t^{-1}.
\end{align*}
The proof is complete.
\hfill$\Box$

\begin{lemma}\label{lemma5}
For any $y>0$, we have
\begin{equation}\label{limit-P}
\lim_{t\to\infty}\widetilde{\mathbb{P}}^{-y}[E_t]=1
\end{equation}
 and
\begin{equation}\label{limit-as}
\lim_{t\to\infty}\inf_{k_t^{1/3}\leq u \leq k_t}\widetilde{\mathbb{P}}^{-y}[E_t \vert \xi^{\lambda_0}_{k_t} = u]  = 1.
\end{equation}
\end{lemma}
\textbf{Proof:}
First, by Lemma \ref{lemma3},
\begin{equation}
\lim_{t\to\infty}\sup_{u\in[k_t^{1/3}, k_t]}\widetilde{\mathbb{P}}^{-y}\left[E_{t,2}^c\Big\vert \xi^{\lambda_0}_{k_t}=u \right]= 0. \label{conclu-1}
\end{equation}
By Lemma \ref{lemma4}, we have
$$
\lim_{t\to\infty}\sup_{u \in [k_t^{1/3}, k_t]}\widetilde{\mathbb{P}}^{-y}\left[E_{t,1}\cap E_{t,2} \cap E_{t,3}^c \big\vert \xi^{\lambda_0}_{k_t} =u \right] = 0.
$$
Note that
\begin{equation}\label{decom-Omega}
\Omega = E_t \cup E_{t,2}^c \cup E_{t,1}^c \cup (E_{t,1}\cap E_{t,2} \cap E_{t,3}^c).
\end{equation}
To prove \eqref{limit-as}, we only need to prove that
\begin{equation}
\inf_{u \in [k_t^{1/3}, k_t]} \widetilde{\mathbb{P}}^{-y}[E_{t,1} \vert \xi^{\lambda_0}_{k_t} =u] \to 1,\quad\mbox{ as } t\to\infty. \label{conclu1}
\end{equation}
Recall that under $\widetilde{\mathbb{P}}^{-y}, y+\xi^{\lambda_0}_t$ is a Bessel-3 process
starting from $y$. Now let $\eta_t:= \xi^{\lambda_0}_t + y$. Then $(\eta, \widetilde{\mathbb{P}}^{-y})$ is equal in law with $(\eta, \widetilde{\Pi}_{y})$. For any $u\in  [k_t^{1/3}, k_t]$,
by the Markov property and Lemma \ref{lemm8}, we have
\begin{align}\label{equal_8}
\widetilde{\mathbb{P}}^{-y}[E_{t,1} \vert \xi^{\lambda_0}_{k_t} =u] &
\ge \widetilde{\Pi}_{y + u} \left(\min_{r \in [0, t- k_t]} \eta_r \geq k_t^{1/6} + y \right) \nonumber\\
& = \frac{1}{y +u} \Pi_0 \left[\left(B_{t-k_t} + y+u\right)1_{\{\min_{r \in [0, t -k_t]} B_r \geq k_t^{1/6} -u\}} \right].
\end{align}
Set $a = u - k_t^{1/6}\geq 0$. Then
using the fact that $\Pi_0 B_{t \land \tau_{-a }} = 0$ for any $t\geq0$, we have
\begin{equation*}
0 = \Pi_0  B_{(t - k_t) \land \tau_{-a }} = -a \Pi_0(\tau_{-a} < t -k_t ) + \Pi_0 (B_{t - k_t} 1_{\{\tau_{-a} \geq t -k_t \}}).
\end{equation*}
Also note that by Lemma \ref{lemma9},
$$
\Pi_{0}(\tau_{-a}\le t -k_t)  = 2 \int^\infty_{a/\sqrt{t-k_t}} \frac{1}{\sqrt{2\pi}}e^{-x^2/2} \textup{d}x.
$$
Then the right-hand of \eqref{equal_8} is equal to
\begin{align}\label{equal_9}
&\frac{1}{y +u}\Pi_0\left[B_{t-k_t} 1_{\{\tau_{-a} \geq t - k_t\}} + (y+u) 1_{\{\tau_{-a} \geq t - k_t\}}\right]\nonumber\\
  &= 1 - \frac{2(y + k_t^{1/6})}{y + u} \int_{(u - k_t^{1/6})/\sqrt{t - k_t}}^ \infty \frac{1}{\sqrt{2\pi}}e^{-x^2/2} \textup{d} x.
\end{align}
By \eqref{equal_8} and \eqref{equal_9}, we get
$$
\widetilde{\mathbb{P}}^{-y}[E_{t,1} \vert \xi^{\lambda_0}_{k_t} =u] \geq 1 - \frac{2(y + k_t^{1/6})}{y + k_t^{1/3}} \int_{0}^ \infty \frac{1}{\sqrt{2\pi}}e^{-x^2/2} \textup{d} x.
$$
By the assumption on $k_t$, we get \eqref{conclu1}.

Now we  prove \eqref{limit-P}.  We claim that
\begin{equation}\label{toprovelimit}
\widetilde{\mathbb{P}}^{-y} [k_t^{1/3} \leq \xi^{\lambda_0}_{k_t} \leq k_t] = \widetilde{\Pi}_y [k_t^{1/3} + y \leq \eta_{k_t} \leq k_t +y ] \to 1,\quad\mbox{ as } t \to \infty.
\end{equation}
In fact, by Theorem 3.2 of \cite{ShWa},  $\lim_{t \to \infty} \log(\eta_t )/ \log t = 1/2,\  \widetilde{\Pi}_y$-a.s. Using the fact that $k_t \to \infty $ as $t \to \infty$, we get \eqref{toprovelimit} holds. Combining \eqref{toprovelimit} and \eqref{conclu-1}, we have
\begin{equation}\label{limit-E2}
\lim_{t\to\infty}\widetilde{\mathbb{P}}^{-y}[E_{t,2}^c] = 0.
\end{equation}
Combining \eqref{toprovelimit} and \eqref{conclu1}, we have
\begin{equation}\label{limit-E1}
\lim_{t\to\infty}\widetilde{\mathbb{P}}^{-y}[E_{t,1}] = 1.
\end{equation}
It follows from Lemma \ref{lemma4} that
\begin{equation}\label{limit-E3}
\lim_{t\to\infty}\widetilde{\mathbb{P}}^{-y}\left[E_{t,1}\cap E_{t,2} \cap E_{t,3}^c \right] =0.
\end{equation}
Using \eqref{decom-Omega}, and combining \eqref{limit-E2}-\eqref{limit-E3}, we obtain \eqref{limit-P}.
\hfill$\Box$

\begin{lemma}\label{lemma11}
For any $y>0$, it holds that
$$
\limsup_{t\to\infty } t \widetilde{\mathbb{P}}^{-y}\left[\frac{ \widetilde{W}_t^{-y}}{\widetilde{W}_t^{-y} +\widetilde{V}_t^{-y}} \frac{1_{E_t}}{\xi^{\lambda_0}_t + y} \right] \leq 	\frac2\pi.
$$
\end{lemma}
\textbf{Proof:}
First note that
$$
\widetilde{\mathbb{P}}^{-y}\left[\frac{\widetilde{W}_t^{-y}}{\widetilde{W}_t^{-y} + \widetilde{V}_t^{-y}} \frac{1_{E_t}}{\xi_t^{\lambda_0} + y} \right]=\widetilde{\mathbb{P}}^{-y}\left[\frac{\widetilde{W}_t^{-y,[k_t, t]}}{\widetilde{W}_t^{-y} + \widetilde{V}_t^{-y}} \frac{1_{E_t}}{\xi^{\lambda_0}_t + y} \right]
+\widetilde{\mathbb{P}}^{-y}\left[\frac{\widetilde{W}_t^{-y,[0, k_t)}}{\widetilde{W}_t^{-y} + \widetilde{V}_t^{-y}} \frac{1_{E_t}}{\xi^{\lambda_0}_t + y} \right].$$
For the first term on the right hand, we have
\begin{equation*}
\widetilde{\mathbb{P}}^{-y}\left[\frac{\widetilde{W}_t^{-y,[k_t, t]}}{\widetilde{W}_t^{-y} + \widetilde{V}_t^{-y}} \frac{1_{E_t}}{\xi^{\lambda_0}_t + y} \right] \leq \widetilde{\mathbb{P}}^{-y}\left[\frac{1/t^2}{\widetilde{V}_t^{-y}(y + k_t^{1/6})} \right] = \frac{1}{y t^2 (k_t^{1/6} +y)},
\end{equation*}
here we used the property that $E_t \subset \{\xi_{t} \geq k_t^{1/6}\},\ E_t \subset E_{t,3}$
 and the equality $\widetilde{\mathbb{P}}^{-y}\left[\frac{1}{\widetilde{V}_t^{-y}}\right]={\mathbb{Q}}^{-y}\left[\frac{1}{{V}_t^{-y}}\right]
=\frac{1}{y}$.
Hence,
$$
\lim_{t\to \infty}t\widetilde{\mathbb{P}}^{-y}\left[\frac{\widetilde{W}_t^{-y,[k_t, t]}}{\widetilde{W}_t^{-y} + \widetilde{V}_t^{-y}} \frac{1_{E_t}}{\xi^{\lambda_0}_t + y} \right]=0.$$
Therefore, we only need to prove that
\begin{equation}\label{toprove-limitsup}
\limsup_{t\to\infty } t \widetilde{\mathbb{P}}^{-y}\left[\frac{ \widetilde{W}_t^{-y, [0,k_t)}}{\widetilde{W}_t^{-y} +\widetilde{V}_t^{-y}} \frac{1_{E_t}}{\xi^{\lambda_0}_t + y} \right] \leq \frac2\pi.
\end{equation}
Note that
\begin{align}\label{ineq4_8}
		&	\widetilde{\mathbb{P}}^{-y}\left[ \frac{\widetilde{W}_t^{-y,[0,k_t)}}{\widetilde{W}_t^{-y} +\widetilde{V}_t^{-y}}  \frac{1_{E_t}}{\xi^{\lambda_0}_t + y} \right]   \leq \widetilde{\mathbb{P}}^{-y}\left[\frac{\widetilde{W}_t^{-y,[0,k_t)}}{\widetilde{W}_t^{-y,[0,k_t)} +\widetilde{V}_t^{-y,[0, k_t)}} \frac{1_{E_t}}{\xi^{\lambda_0}_t + y} \right] \nonumber\\
&\leq  \widetilde{\mathbb{P}}^{-y}\left[\frac{\widetilde{W}_t^{-y,[0,k_t)}}{\widetilde{W}_t^{-y,[0,k_t)} +\widetilde{V}_t^{-y,[0, k_t)}} 1_{\{\xi^{\lambda_0}_{k_t} \in [k_t^{1/3}, k_t] \}} \right] \times \sup_{u \in [k_t^{1/3}, k_t]} \widetilde{\mathbb{P}}^{-y} \left[\frac{1}{\xi^{\lambda_0}_t +y }\bigg\vert \xi^{\lambda_0}_{k_t} = u\right].
\end{align}
In the last inequality we used the Markov property of $\xi$.  Let $\{(\eta_t)_{t \geq 0}, \widetilde{\Pi}_{u+y}\}$ be a Bessel-3 process starting from $u+y$.  By Lemmas \ref{lemm8} and \ref{lemma9}, we have
\begin{align}\label{ineq4_9}
&\widetilde{\mathbb{P}}^{-y} \left[\frac{1}{\xi^{\lambda_0}_t +y }\bigg\vert \xi^{\lambda_0}_{k_t} = u\right]  = \widetilde{\Pi}_{u+y}\left[\frac{1}{\eta_{t-k_t}} \right] = \frac{1}{u+y}\Pi_{u+y}\left[1_{\{\min_{r \in [0, t -k_t] } B_r > 0 \}}
\right] \nonumber\\
&= \frac{1}{u + y} \Pi_0( \tau_{-(y +u)} > t -k_t)
 = \frac{2}{y+u}\int_0^{(y+u)/\sqrt{t -k_t}} \frac{1}{\sqrt{2\pi}}e^{-x^2/2} \textup{d}x .
\end{align}
By \eqref{ineq4_8} and \eqref{ineq4_9}, we get
\begin{align}\label{ineq4_28}
\widetilde{\mathbb{P}}^{-y}\left[ \frac{\widetilde{W}_t^{-y,[0,k_t)}}{\widetilde{W}_t^{-y} +\widetilde{V}_t^{-y}}   \frac{1_{E_t}}{\xi^{\lambda_0}_t + y} \right]  \leq & \widetilde{\mathbb{P}}^{-y}\left[\frac{\widetilde{W}_t^{-y,[0,k_t)}}{\widetilde{W}_t^{-y,[0,k_t)} + \widetilde{V}_t^{-y,[0, k_t)}} 1_{\{\xi^{\lambda_0}_{k_t} \in [k_t^{1/3}, k_t] \}} \right] \nonumber\\
& \times  \sup_{u \in [k_t^{1/3}, k_t]} \frac{2}{y+u}\int_0^{(y+u)/\sqrt{t -k_t}} \frac{1}{\sqrt{2\pi}}e^{-x^2/2} \textup{d}x .
\end{align}
Because $\lim_{\varepsilon \to 0^+} \frac{2}{\varepsilon} \int_0^\varepsilon e^{-x^2/2}/\sqrt{2\pi}\textup{d}x = \sqrt{2/\pi}$	
and $(y+u)/\sqrt{t-k_t}$ converges to $0$ uniformly on $u\in [k_t^{1/3}, k_t]$ as $t \to \infty $, we have
\begin{equation}
\sup_{u \in [k_t^{1/3}, k_t]} \frac{2\sqrt{t}}{y+u}\int_0^{(y+u)/\sqrt{t -k_t}} \frac{1}{\sqrt{2\pi}}e^{-x^2/2} \textup{d}x \to
\sqrt{\frac2\pi}. \label{equal_2}
\end{equation}
Using the Markov property at time $k_t$ again, we get
\begin{align}
&\widetilde{\mathbb{P}}^{-y}\left[\frac{\widetilde{W}_t^{-y,[0,k_t)}}{\widetilde{W}_t^{-y,[0,k_t)} + \widetilde{V}_t^{-y,[0, k_t)}}  1_{E_t}\right]  \nonumber\\
& \geq \widetilde{\mathbb{P}}^{-y}\left[\frac{\widetilde{W}_t^{-y,[0,k_t)}}{\widetilde{W}_t^{-y,[0,k_t)} + \widetilde{V}_t^{-y,[0, k_t)}} 1_{\{\xi^{\lambda_0}_{k_t} \in [k_t^{1/3}, k_t] \}} \right] \cdot \inf_{u \in [k_t^{1/3},k_t]} \widetilde{\mathbb{P}}^{-y}[E_t \vert \xi^{\lambda_0}_{k_t} = u].\label{ineq4_11}
\end{align}
	Because $\widetilde{W}_t^{-y,[0,k_t)}/(\widetilde{W}_t^{-y,[0,k_t)} + \widetilde{V}_t^{-y,[0, k_t)}) \cdot 1_{E_t} \leq 1$,
the left-hand of \eqref{ineq4_11} is bounded  above by
\begin{align}
		&\widetilde{\mathbb{P}}^{-y}\left[\frac{\widetilde{W}_t^{-y,[0,k_t)}}{\widetilde{W}_t^{-y,[0,k_t)} +\widetilde{V}_t^{-y,[0, k_t)}}  1_{E_t}\right]	
\leq  \widetilde{\mathbb{P}}^{-y}\left[\frac{\widetilde{W}_t^{-y,[0,k_t)}}{\widetilde{W}_t^{-y,[0,k_t)} +\widetilde{V}_t^{-y,[0, k_t)}}  1_{E_t}1_{\{\widetilde{V}_t^{-y} > 1/t\}}\right] + \widetilde{\mathbb{P}}^{-y}\left[\widetilde{V}_t^{-y} \leq \frac{1}{t}\right] \nonumber\\
		& \leq \widetilde{\mathbb{P}}^{-y}\left[\frac{\widetilde{W}_t^{-y,[0,k_t)}}{\widetilde{V}_t^{-y,[0, k_t)}}  1_{E_t}1_{\{\widetilde{V}_t^{-y} > 1/t\}}\right]+  \frac{1}{t} \widetilde{\mathbb{P}}^{-y}\left[\frac{1}{\widetilde{V}_t^{-y}}\right]
		= \widetilde{\mathbb{P}}^{-y}\left[\frac{\widetilde{W}_t^{-y,[0,k_t)}}{\widetilde{V}_t^{-y,[0, k_t)}}  1_{E_t}1_{\{\widetilde{V}_t^{-y} > 1/t\}}\right]+ \frac{1}{ty}, \label{ineq4_12}
\end{align}
where in the last inequality we used the Markov inequality for $\left(\widetilde{V}_t^{-y} \right)^{-1}.$
Fix a constant $\eta \in (0,1)$,
on $E_t \cap \{ \widetilde{V}_t^{-y} > 1/t\}$,
we have, for large $t$
such that $t > \eta^{-1}$,
$\widetilde{V}_t^{-y,[k_t, t]} \leq \eta \widetilde{V}_t^{-y}.$
So when $t$ is large, using \eqref{ineq4_12}, we have
$$
\widetilde{\mathbb{P}}^{-y} \left[\frac{\widetilde{W}_t^{-y,[0,k_t)}} {\widetilde{W}_t^{-y,[0,k_t)} + \widetilde{V}_t^{-y,[0, k_t)}}  1_{E_t}\right] \leq \frac{1}{ty} + \frac{1}{1-\eta} \widetilde{\mathbb{P}}^{-y} \left[\frac{\widetilde{W}^{-y}_t}{\widetilde{V}_t^{-y}} \right].
$$
By \eqref{limit-W-V}, we have
\begin{equation}\label{ineq4_13}
\widetilde{\mathbb{P}}^{-y} \left[\frac{\widetilde{W}_t^{-y,[0,k_t)}} {\widetilde{W}_t^{-y,[0,k_t)} + \widetilde{V}_t^{-y,[0, k_t)}}  1_{E_t}\right]
	\leq
	\frac{\sqrt{2/\pi}}{(1-\eta)\sqrt{t}} + o\left(\frac{1}{\sqrt{t}} \right), \quad\mbox{ as } t\to \infty.
\end{equation}
By \eqref{limit-as}, \eqref{ineq4_28}, \eqref{equal_2}, \eqref{ineq4_11} and \eqref{ineq4_13}, we finally get that
$$
\limsup_{t\to\infty } t \widetilde{\mathbb{P}}^{-y}\left[\frac{ \widetilde{W}_t^{-y,[0,k_t) }}{\widetilde{W}_t^{-y} + \widetilde{V}_t^{-y}} \frac{1_{E_t}}{\xi^{\lambda_0}_t + y} \right]
\leq \frac{2}{\pi(1-\eta)}.
$$
Since  the above holds for any small $\eta\in(0,1)$, \eqref{toprove-limitsup} holds. The proof is complete.
\hfill$\Box$

\textbf{Proof of Proposition \ref{convergence-L2}:}
Applying Lemmas \ref{lemma10} and \ref{lemma11}, and \eqref{step_6}, we get
\begin{align*}
		&	\limsup_{t\to \infty} \widetilde{\mathbb{P}}^{-y}\left[\left(\frac{\sqrt{t}\widetilde{W}_t^{-y}}{\widetilde{V}_t^{-y}+\widetilde{W}_t^{-y}} - \sqrt{\frac2\pi}  \right)^2 \right] \\ &=\limsup_{t\to \infty} \left\{ \widetilde{\mathbb{P}}^{-y}\left[\left(\frac{\sqrt{t}\widetilde{W}_t^{-y}}{\widetilde{V}_t^{-y}+\widetilde{W}_t^{-y}} \right)^2\right] - \frac2\pi \right\} - 2\sqrt{\frac2\pi} \lim_{t \to \infty}\left\{ \widetilde{\mathbb{P}}^{-y}\left[\frac{\sqrt{t}\widetilde{W}_t^{-y}}{\widetilde{V}_t^{-y}+\widetilde{W}_t^{-y}}\right] - \sqrt{\frac2\pi} \right\}  \leq 0,
\end{align*}
which means that \eqref{ineq4_22} holds.
\hfill$\Box$

\textbf{Proof of Theorem \ref{th1}:}
Let  $\mathcal{R}^{\lambda_0}$ and $\widetilde{\mathcal{R}}^{\lambda_0}$ be the smallest closed set containing  $\bigcup_{t\ge 0}\mbox{supp}X_t^{\lambda_0}$  and $\bigcup_{t\ge 0}\mbox{supp}\widetilde X_t^{\lambda_0}$, respectively.
Then by \cite[ Corollary 3.2]{KLMR}, under condition \eqref{condition1},
$\mathbb{P}(\inf \mathcal{R}^{\lambda_0} > -\infty ) = 1.$
So for any $0 < \eta < \mathbb{P}(\mathcal{E}^c)$, there exists $K>0$ such that $\mathbb{P}(\inf \mathcal{R}^{\lambda_0} > - K) > 1 -\eta.$
Let $y: = K $ be fixed and define $\Omega_k:=\{\inf \mathcal{R}^{\lambda_0} > -K\}$ and $\widetilde\Omega_k:=\{\inf \widetilde{\mathcal{R}}^{\lambda_0} > -K\}$. Then
$$
\mathbb{P}(\Omega_K\cap \mathcal{E}^c) \geq \mathbb{P}(\Omega_K) + \mathbb{P}(\mathcal{E}^c) -1 > 1-\eta + \mathbb{P}(\mathcal{E}^c)-1 > 0.
$$
For any $\varepsilon> 0$, put
$$
G_t = \left\{\Big\vert \frac{\sqrt{t} W_t^{-y}}{V_t^{-y}+W_t^{-y}} - \sqrt{\frac2\pi}    \Big\vert > \varepsilon \right\},\quad
 \widetilde{G}_t = \left\{\Big\vert \frac{\sqrt{t} \widetilde{W}_t^{-y}}{\widetilde{V}_t^{-y}+\widetilde{W}_t^{-y}} - \sqrt{\frac2\pi}    \Big\vert > \varepsilon \right\}.
$$
Define $\mathbb{P}^{**} (\cdot) = \mathbb{P}(\cdot \vert \Omega_K\cap \mathcal{E}^c)$.
By \eqref{ineq4_22} we have $\lim_{t \to \infty} \widetilde{\mathbb{P}}^{-y}[\widetilde{G}_t] =0$. Thus,
$$
\frac{\mathbb{P}(\Omega_K\cap \mathcal{E}^c)}{y}\lim_{t \to \infty} \mathbb{P}^{**}[V_t^{-y}1_{G_t}] =\lim_{t\to\infty} \widetilde{\mathbb{P}}^{-y}[\widetilde{G}_t\cap \widetilde{\Omega}_K \cap \widetilde{\mathcal{E}}^c] = \lim_{t \to \infty}\widetilde{\mathbb{P}}^{-y}[\widetilde{G}_t] = 0,
$$
where $\widetilde{\cal E}:=\{\exists t\geq 0\mbox{ such that} \|\widetilde X_t^{\lambda_0}\|=0\}$ with  $\widetilde{\mathbb{P}}^{-y}$-probability $0$.
Then by Proposition \ref{convergence-L2}, we have
\begin{equation}\label{equal_10}
V_t^{-y}1_{G_t} \xrightarrow[t \to \infty]{\mathbb{P}^{**}} 0.
\end{equation}
Notice that on the event $\Omega_K := \{\inf \mathcal{R}^{\lambda_0} > - K \}$,
we have
$$
V_t^{-y} = V_t^{-K} = \partial W_t + K W_t > 0,\ \ \ \ W_t^{-y} = W_t^{-K} = W_t,$$
and $\lim_{t\to\infty} V_t^{-y} = \partial W_\infty>0\ \mathbb{P}^{**}$-a.s..
Together with \eqref{equal_10} we get $\lim_{t\to\infty}\mathbb{P}^{**}[G_t] = 0$ for any $\varepsilon > 0$, which says
\begin{equation}
\frac{\sqrt{t} W_t^{-y}}{V_t^{-y}+W_t^{-y}} =	\frac{\sqrt{t}W_t}{\partial W_t + (K+1) W_t}\xrightarrow[t \to \infty]{\mathbb{P}^{**}} \sqrt{\frac2\pi}. \label{equal_3}
\end{equation}
Recall that $\mathbb{P}(\mathcal{E}^c) = 1 - e^{-\lambda^*} >0 $ and
$\mathbb{P}^{**}(W_t > 0, \,\forall t>0)=\mathbb{P}^{**}(\lim_{t \to \infty} W_t > 0)=1$.
According to \eqref{equal_3} we get
$$	
\frac{\partial W_t}{\sqrt{t}W_t}\xrightarrow[t \to \infty]{\mathbb{P}^{**}} \sqrt{\frac\pi2}.
$$
For any $\gamma >0$, define
$$
A_t = \left\{\Big\vert \frac{\partial W_t}{\sqrt{t}W_t }  - \sqrt{\frac\pi2} \Big\vert > \gamma \right\}.
$$
Then
$\lim_{t \to \infty}\mathbb{P}^{**}[1_{A_t}]=0.$
Noticing that $\mathbb{P}^{*} (\cdot) = \mathbb{P}(\cdot \vert \mathcal{E}^c)$
and
$\mathbb{P}^{*}[1_{A_t}1_{\Omega_K}] = \mathbb{P}^{**}[1_{A_t}] \mathbb{P}(\Omega_K\cap \mathcal{E}^c)/\mathbb{P}(\mathcal{E}^c)$, we obtain that
$$
1_{A_t}1_{\Omega_K} \xrightarrow[t \to \infty]{\mathbb{P}^{*}} 0,
$$
which means	$\limsup_{t\to\infty } \mathbb{P}^{*}(A_t) \leq \lim_{t \to \infty} \mathbb{P}^*(A_t\cap \Omega_K) + \mathbb{P}^*(\Omega_K^c) \leq \eta/\mathbb{P}(\mathcal{E}^c). $	Since $\eta$ is arbitrary, we deduce that $\lim_{t \to \infty} \mathbb{P}^*(A_t)=0$ for any $\gamma > 0$,
 which says
$$
\frac{\partial W_t}{\sqrt{t}W_t }  \xrightarrow[t \to \infty]{\mathbb{P}^{*}} \sqrt{\frac\pi2}.
$$
This is also equivalent to say that,  on the event $\mathcal{E}^c$, we have
\begin{equation}
\sqrt{t}W_t  \xrightarrow[t \to \infty]{\mathbb{P}}
\sqrt{\frac{2}{\pi}}\partial W_\infty .\label{final_3}
\end{equation}
On $\mathcal{E}$, \eqref{final_3} holds obviously, and the proof is now complete.
\hfill$\Box$

\section{Proof of Theorem \ref{th2}}
Recall the definitions of the process $\{(Z_t, \Lambda_t)_{t\geq 0}\}$ and the probability measures $\mathbf{P}_{(\mu, \eta)}$ and   $\mathbf{P}_\mu$
 with $\mu\in \mathcal{M}(\R)$ and $\eta\in\mathcal{M}_a(\R)$, defined in Subsection \ref{Skeleton}. Set $\mathbf{P}:= \mathbf{P}_{\delta_0}$.
By the skeleton decomposition for $X$,
$(\Lambda_t, \mathbf{P})$ is equal in law to $(X,\P)$. To prove Theorem \ref{th2}, we only need to prove that on survival event $\left(\mathcal{E}^{\Lambda}\right)^c$ where $\mathcal{E}^\Lambda : = \{\lim_{t\to\infty}\Vert \Lambda_t \Vert =0 \}$,
 \begin{equation}\label{liminf-Lambda}
 \limsup_{t\to\infty} \sqrt{t}\langle e^{-\lambda_0 (\cdot +\lambda_0 t)}, \Lambda_t\rangle= +\infty \quad\mathbf{P}\mbox{-almost surely}.
 \end{equation}
The intuitive idea for proving the limit above is that the behaviour of $\Lambda$ is determined by the skeleton $Z$. By branching property of $Z$ we only consider the law $\mathbf{P}_{(\delta_0,\delta_0)}.$ Let $\{\mathbf{e}_n: n\ge 1\}$ be iid exponential random variables independent of $Z$. Let $T_0:=0$ and
$T_n=\sum^n_{i=1}\mathbf{e}_i$ for $n\ge 1$.
If we look at $Z$ at independent times  $\{T_n:n=1,2,...\}$, then  $\{Z_{T_n}, n\geq 1\}$ is a branching random walk.
We expect the  behavior of this branching random walk to dominate the behavior of $\Lambda$. Let $\{\mathcal{Z}_{n}, n\geq 1\}$  be the translation of $\{Z_{T_n}, n\geq 1\}$ defined in \eqref{def-mathcal-Z} below. We will show that $\{\mathcal{Z}_{n}, n\geq 1\}$ satisfies conditions of  Aidekon and  Shi \cite{AESZ}. Then by \cite[Theorem 6.1]{AESZ},
$$
\liminf_{n\to\infty} \left(L^{\mathcal Z}_n-\frac{1}{2}\log n \right)= -\infty\quad\mathbf{P}_{(\delta_0, \delta_0 )}\mbox{-almost surely},
$$
where $L^{\mathcal Z}_n$ is minimum of the support  of $\mathcal{Z}_n$.   Let   $L_t^{Z}$ be  minimum of the support of $Z_t$.
By  definition  \eqref{def-mathcal-Z},
$L^{\mathcal Z}_n=\lambda_0 (L_{T_n}^{Z} +\lambda_0 T_n)$,
and then we have
\begin{equation}\label{liminf-LZ}
\liminf_{n\to\infty} \left(\lambda_0 (L_{T_n}^{Z}
+\lambda_0 T_n)
-\frac{1}{2}\log T_n \right)= -\infty\quad\mathbf{P}_{(\delta_0,\delta_0)}\mbox{-almost surely}.
\end{equation}
We will bound $\langle e^{-\lambda_0 (\cdot +\lambda_0 T_n)}, \Lambda_{T_n}\rangle$ from below by immigrations along the path of  $L_{\cdot}^{Z}$, and then
use the limit result \eqref{liminf-LZ} for $L_{T_n}^{Z}$ to get \eqref{liminf-Lambda}.

Now we give a more precise proof.
Note that
\begin{equation}\label{Sum1}
\mathbf{P}(\cdot)= \sum_{k=0}^\infty \frac{(\lambda^*)^{k}}{k !}e^{-\lambda^*} \mathbf{P}_{(\delta_0, k\delta_0)}(\cdot),
\end{equation}
and $\mathbf{P}(\mathcal{E}^\Lambda) = \P(\mathcal{E})=e^{-\lambda^*}$.
It is obvious that $\mathbf{P}_{(\delta_0, 0\delta_0)}(\mathcal{E}^\Lambda)=1$.
Together with \eqref{Sum1}, we know that for $k\geq 1, \mathbf{P}_{(\delta_0, k\delta_0)}(\mathcal{E}^\Lambda) = 0$.
Thus, to prove Theorem \ref{th2}, it suffices to show that, for any $k\ge 1$,
the limsup in \eqref{e:th2} is valid $\mathbf{P}_{(\delta_0, k\delta_0)}$-almost surely.
By the branching property, without loss of generality, we only need to deal with the case of $k=1$.

Let $\{\mathbf{e}_n: n\ge 1\}$ be iid exponential random variables with parameter $\kappa\in (0,\infty)$, independent of $Z$. Let $T_0:=0$ and
$T_n=\sum^n_{i=1}\mathbf{e}_i$ for $n\ge 1$.
Now for $n\geq 1$, we define $\mathcal{Z}_n$ so that, for any $f\in \B_b^+(\R)$,
\begin{equation}\label{def-mathcal-Z}
\langle f, \mathcal{Z}_n \rangle = \langle f\left(\lambda_0(\cdot +\lambda_0 T_n)\right),Z_{T_n}\rangle.
\end{equation}
Then  $\{(\mathcal{Z}_n )_{n\geq 1}, \mathbf{P}_{(\delta_0, \delta_0)}\}$ is a branching random walk.
By \eqref{Generating-F},
define $m:=\sum_{n\geq0} np_n = F'(1-)$.
It is easy to check that $\lambda_0 =\sqrt{2\psi'(\lambda^*)(m-1) }$. We first check that the conditions of \cite[Theorem 6.1]{AESZ} for $\mathcal{Z}$ are satisfied. More precisely, under assumption \eqref{condition2},
\eqref{moment1} \eqref{moment2} and \eqref{moment3} hold.
For simplicity, we define
$$
W_n^\mathcal{Z}:= \langle e^{-\cdot}, \mathcal{Z}_n \rangle,\quad D_n^\mathcal{Z}:= \langle \cdot e^{-\cdot}, \mathcal{Z}_n \rangle,\quad  D_n^{\mathcal{Z},2}:= \langle (\cdot)^2 e^{-\cdot}, \mathcal{Z}_n \rangle,\quad D_n^{\mathcal{Z},+}:= \langle (\cdot)_+ e^{-\cdot}, \mathcal{Z}_n \rangle.
$$
The additive martingale associated to $Z$ with parameter $\lambda$ is defined as
\begin{equation}\label{def-add-mart}
W_s^{Z}(\lambda):= e^{-\lambda c_\lambda s}\langle e^{-\lambda \cdot}, Z_s\rangle = e^{-(\lambda - \lambda_0)^2 s/2}\langle e^{-\lambda(\cdot +\lambda_0 s)}, Z_s\rangle,
\end{equation}
where $c_\lambda:= \lambda/2 + \psi'(\lambda^*)(m-1)/\lambda = (\lambda^2 + \lambda_0^2)/(2\lambda)$ and $\lambda c_\lambda = (\lambda - \lambda_0)^2 /2 + \lambda \lambda_0.$

\begin{lemma}\label{Integrability-condition1}
If  $\sum_{n\geq 1} n(\log n)^2 p_n < \infty$, then
\begin{equation}\label{Check1}
\mathbf{P}_{(\delta_0, \delta_0)}\left[W_1^{\mathcal{Z}}\right]=1,\quad 	\mathbf{P}_{(\delta_0, \delta_0)}\left[D_1^{\mathcal{Z}}\right]=0, \quad
\mathbf{P}_{(\delta_0, \delta_0)}\left[D_1^{\mathcal{Z},2}\right] < \infty
\end{equation}
and
\begin{equation}\label{Check2}
\mathbf{P}_{(\delta_0, \delta_0)}\left[W_1^{\mathcal{Z}}\log_+^2 W_1^{\mathcal{Z}}\right]<\infty,\quad \mathbf{P}_{(\delta_0, \delta_0)}\left[D_1^{\mathcal{Z},+}\log_+ D_1^{\mathcal{Z},+}\right]<\infty.
\end{equation}
\end{lemma}
\textbf{Proof :} $Step\ 1:$ Define $W_s^{Z}$ and $ D_s^{Z}$ by
$$
W_s^{Z}:= \langle e^{-\lambda_0(\cdot +\lambda_0 s)}, Z_s\rangle, \quad  D_s^{Z}:= \langle (\cdot + \lambda_0 s) e^{-\lambda_0(\cdot +\lambda_0 s)}, Z_s\rangle.
$$
Then by \cite{Ky}, $W_s^{Z}$ and $ D_s^{Z}$ are the additive martingale and the derivative martingale associated to the branching Brownian motion $Z$ in the critical case $\lambda =\lambda_0$ respectively.
		
By some direct calculation and the martingale property, we have
\begin{align*}
\mathbf{P}_{(\delta_0, \delta_0)}\left[W_1^{\mathcal{Z}}\right] &= \int_0^\infty \kappa e^{-\kappa s}\mathbf{P}_{(\delta_0, \delta_0)}\left[W_s^{Z}\right]\mathrm{d} s = \int_0^\infty \kappa e^{-\kappa s}\mathrm{d}s =1,\\
\mathbf{P}_{(\delta_0, \delta_0)}\left[D_1^{\mathcal{Z}}\right] & = \int_0^\infty \kappa e^{-\kappa s}\mathbf{P}_{(\delta_0, \delta_0)}\left[D_s^{Z}\right]\mathrm{d}s =0.
\end{align*}
Now define
$$
D_s^{Z,2}:= \lambda_0^2\langle (\cdot+\lambda_0 s)^2 e^{-\lambda_0(\cdot+\lambda_0 s)}, Z_s \rangle.
$$
Using the many-to-one formula, we get
\begin{align*}
&\mathbf{P}_{(\delta_0, \delta_0)}\left[D_1^{\mathcal{Z},2}\right] = \int_0^\infty \kappa e^{-\kappa s}\mathbf{P}_{(\delta_0, \delta_0)}\left[D_s^{Z,2}\right]\mathrm{d}s
= \int_0^\infty \kappa e^{-\kappa s} \lambda_0^2e^{\lambda_0^2s /2}\Pi_0 \left[(B_s+\lambda_0 s)^2 e^{-\lambda_0(B_s + \lambda_0s)}\right]\mathrm{d}s \\  &=\lambda_0^2 \int_0^\infty \kappa e^{-\kappa s}  \Pi_0^{-\lambda_0} \left[(B_s+\lambda_0 s)^2 \right]\mathrm{d}s= \lambda_0^2 \int_0^\infty \kappa s e^{-\kappa s} \mathrm{d}s < \infty.
\end{align*}
Thus, \eqref{Check1} holds.
		
$Step\ 2:$
In this step we prove the first inequality of \eqref{Check2}.
Define a new probability $\Q^{Z}$ by
$$
\frac{\mathrm{d} \Q^{Z}}{\mathrm{d}\mathbf{P}_{(\delta_0, \delta_0)}}\bigg|_{\sigma(Z_r^{1}, r\leq s)} := W_s^{Z},\quad s\geq 0.
$$
Then under $\Q^{Z}, Z$ has the following spine decomposition:
		
(i) There is a initial marked particle moving as a Brownian motion with drift $-\lambda_0$ starting from $0$, we denote the trajectory of this particle by $w_s$.
		
(ii) The branching rate of this marked particle is $\psi'(\lambda^*)m$ and the
offspring distribution of the marked particle is given by $\widetilde{p}_n:= np_n/m,
n=1, 2, \dots$.
		
(iii) When the marked particle dies, given the number of the offspring, mark one of its offspring uniformly.
		
(iv) The unmarked individuals evolve independently as $Z$ under $\P_{(\delta_0, \delta_0)}$.
		
Note that
\begin{equation}\label{Equal-0}
\mathbf{P}_{(\delta_0, \delta_0)}\left[W_1^{\mathcal{Z}}\log_+^2 W_1^{\mathcal{Z}}\right] = \int_0^\infty \kappa e^{-\kappa s}\mathbf{P}_{(\delta_0, \delta_0)}\left[W_s^{Z}\log_+^2 W_s^{Z}\right]\mathrm{d}s.
\end{equation}
By a change of measure, we have
$$
\mathbf{P}_{(\delta_0, \delta_0)}\left[W_s^{Z}\log_+^2 W_s^{Z}\right] = \Q^{Z} \left[\log_+^2 W_s^{Z}\right].
$$
Let $A > 4$ be a constant such that
\begin{equation}\label{Ineq-1}
\log A (\log A - 2\log 2) \geq \sup_{a \geq 1} \left(\log^2 (a+1) - \log^2 a\right),
\end{equation}
There exists such an $A$ since for all $a \geq 1$, by inequality $\ln (x+1)\leq x$, we have
$$
\log_+^2 (a+1) - \log_+^2 a = \left(\log(a+1)+ \log a\right)\left(\log\left(1+a^{-1}\right)\right)\leq (2a-1)\times a^{-1}<2.
$$
Now let $b,c \geq A$, using \eqref{Ineq-1}, it is easy to check that the inequality
\begin{equation}\label{Ineq-2}
\log^2 (b+c) \leq \log^2 b + \log^2 c
\end{equation}
holds by assuming $b\geq c$ and $b =ac$. For $\ell\ge 1$,
we use $\Gamma_\ell$ to denote the $\ell$-th fission time of the spine under
$\Q^{Z}$ and $O_\ell $ the number of offspring at the fission time $\Gamma_\ell$.
Then
\begin{align}\label{Equal-1}
W_s^{Z}= &\sum_{\ell\geq 1} 1_{\{\Gamma_\ell \leq s  \}} e^{-\lambda_0^2 \Gamma_\ell}W_{s-\Gamma_\ell}^{Z,\Gamma_\ell} 1_{\left\{ e^{-\lambda_0^2 \Gamma_\ell} W_{s-\Gamma_\ell}^{Z,\Gamma_\ell}<A\right\}} \nonumber\\
& \quad
+\sum_{\ell\geq 1} 1_{\{\Gamma_\ell \leq s  \}}e^{-\lambda_0^2 \Gamma_\ell} W_{s-\Gamma_\ell}^{Z,\Gamma_\ell} 1_{\left\{ e^{-\lambda_0^2 \Gamma_\ell}W_{s-\Gamma_\ell}^{Z,\Gamma_\ell}\geq A\right\}}+ e^{-\lambda_0(w_s + \lambda_0s)} \nonumber\\
=&: H_1 + H_2 + H_3 ,
\end{align}
where, given the information along the spine, $W^{Z,\Gamma_\ell}$ is the additive martingale
associated with the branching Brownian motion starting from the $O_\ell-1$ unmarked individuals.
Note that for any  $x,y,z >0$,	we have $\log_+^2 (x+y+z)\leq \log_+^2 (3x) + \log_+^2 (3y) + \log_+^2 (3z)$ and $\log_+^2 x \leq 4 x $. Then \eqref{Equal-1} implies that
\begin{equation}\label{Ineq-3}
\log_+^2 W_s^{Z} \leq \log_+^2 (3H_1) +  \log_+^2 (3H_2) + \log_+^2 (3H_3)  \leq 12 H_1  + \log_+^2 (3H_2) +\log_+^2 (3H_3).
\end{equation}
Since $H_1 \leq A \sum_{\ell\geq 1} 1_{\{\Gamma_\ell \leq s\}}$, we have
\begin{equation}\label{Ineq-4}
\Q^{Z}[H_1] \leq A\int_0^s \psi'(\lambda^*) m\mathrm{d}r = A \psi'(\lambda^*) ms.
\end{equation}
Also, note that $w_s +\lambda_0 s$ under $\Q^{Z}$ is a standard Brownian motion, so
\begin{align}\label{Ineq-5}
&\Q^{Z}\left[\log_+^2 (3H_3)\right] \leq  2(\log 3)^2 + 2\Q^{Z}\left[\log_+^2 (H_3)\right] \nonumber\\
&\leq 2(\log 3)^2  + 2\lambda_0^2 \Q^{Z}(w_s +\lambda_0 s)^2 = 2(\log 3)^2  + 2\lambda_0^2 s.
\end{align}
Here in the first inequality above we used inequality
\begin{equation}\label{int-log}
\log_+^2(ab)\leq (\log_+ a + \log_+ b)^2 \leq 2\log_+^2 a + 2\log_+^2 b.\end{equation}
Define
$$
\overline{W}_{s-\Gamma_\ell}^{Z,\Gamma_\ell} := e^{\lambda_0w_{\Gamma_\ell} } {W}_{s-\Gamma_\ell}^{Z,\Gamma_\ell}.
$$
Using \eqref{Ineq-2} and \eqref{int-log} again,
we deduce that
\begin{align}\label{Ineq-10}
&\log_+^2 (3H_2) \leq  2(\log 3)^2 + 2\log_+^2 (H_2)\nonumber\\	
&\leq 2(\log 3)^2 + 2 \sum_{\ell \geq 1} 1_{\{\Gamma_\ell \leq s\}}1_{\left\{  e^{-\lambda_0^2 \Gamma_\ell}W_{s-\Gamma_\ell}^{Z,\Gamma_\ell}\geq A\right\}} \log_+^2  \left[e^{-\lambda_0^2 \Gamma_\ell}W_{s-\Gamma_\ell}^{Z,\Gamma_\ell}\right] \nonumber\\
&\leq  2(\log 3)^2 + 4\sum_{\ell \geq 1} 1_{\{\Gamma_\ell \leq s\}} \log_+^2 \overline{W}_{s-\Gamma_\ell}^{Z,\Gamma_\ell} + 4\sum_{\ell\geq 1} 1_{\{\Gamma_\ell \leq s\}} \log_+^2 \left(e^{-\lambda_0(w_{\Gamma_\ell} + \lambda_0 \Gamma_\ell)}\right) \nonumber\\
&\leq  2(\log 3)^2 + 4\sum_{\ell \geq 1} 1_{\{\Gamma_\ell \leq s\}} \log_+^2 \overline{W}_{s-\Gamma_\ell}^{Z,\Gamma_\ell} + 4\lambda_0^2\sum_{\ell\geq 1} 1_{\{\Gamma_\ell \leq s\}}(w_{\Gamma_\ell} + \lambda_0 \Gamma_\ell)^2	.
\end{align}
Similarly, we have
\begin{equation}\label{Ineq-6}
\Q^{Z} \left[\sum_{\ell\geq 1} 1_{\{\Gamma_\ell \leq s\}}(w_{\Gamma_\ell} + \lambda_0 \Gamma_\ell)^2\right] = \psi'(\lambda^*)m\int_0^s \Q^{Z}\left[(w_r +\lambda_0r)^2\right]\mathrm{d} r = \psi'(\lambda^*)m s^2/2.
\end{equation}
Now given $w,\Gamma_\ell$ and $O_\ell$, by the spatial homogeneity of branching Brownian motion, we have that $\Q^{Z}\left[\overline{W}_{s-\Gamma_\ell}^{Z,\Gamma_\ell} \big| w,\Gamma_\ell, O_\ell\right] = O_\ell -1.$
By the branching property of $Z$, we have
$\overline{W}_{s-\Gamma_\ell}^{Z,\Gamma_\ell} = \sum_{j=1}^{O_\ell -1}\overline{W}_{s-\Gamma_\ell}^{Z,\Gamma_\ell, j}$, where
$\overline{W}_{s-\Gamma_\ell}^{Z,\Gamma_\ell, j}$, $j=1,\cdots, O_\ell -1$, are independent and have the same distribution given $w, \Gamma_\ell$ and $O_\ell$.
Thus,
\begin{equation}\label{Ineq-7}
\Q^{Z}\left[ \log_+^2  \overline{W}_{s-\Gamma_\ell}^{Z,\Gamma_\ell} \Big|w, \Gamma_\ell, O_\ell  \right]\leq 2\log_+^2(O_\ell -1) + 2 \Q^{Z}\left[  \log_+^2 \left(\max_{j\leq O_\ell -1} \overline{W}_{s-\Gamma_\ell}^{Z,\Gamma_\ell, j} \right)\Big|w, \Gamma_\ell, O_\ell  \right].
\end{equation}
By the Markov inequality,
\begin{align}\label{Ineq-8}
& \Q^{Z}\left[  \log_+^2 \left(\max_{j\leq O_\ell -1} \overline{W}_{s-\Gamma_\ell}^{Z,\Gamma_\ell, j} \right)\Big|w, \Gamma_\ell, O_\ell  \right] = \int_0^\infty 2y \mathrm{d}y \Q^{Z}\left[  \max_{j\leq O_\ell -1} \overline{W}_{s-\Gamma_\ell}^{Z,\Gamma_\ell, j} > e^y\Big|w, \Gamma_\ell, O_\ell  \right]\nonumber\\
&= \int_0^\infty 2y \mathrm{d}y \left[1- \prod_{j\leq O_\ell -1}\left(1- \Q^{Z}\left[   \overline{W}_{s-\Gamma_\ell}^{Z,\Gamma_\ell, j} >  e^y\Big|w, \Gamma_\ell, O_\ell  \right]\right)\right]\nonumber\\
&\leq  \int_0^\infty 2y \mathrm{d}y \left[1- \prod_{j\leq O_\ell -1}(1-e^{-y})\right] = \int_0^\infty 2y \left[1-(1-e^{-y})^{O_\ell -1}\right]\mathrm{d}y.
\end{align}
When $O_\ell -1 < e^{y/2}$,
using the fact that $(1-x)^k \geq 1-kx$  for all $x\le 1$,
we get
$$
2y \left[1-(1-e^{-y})^{O_\ell -1}\right] \leq 2y (O_\ell- 1)e^{-y} \leq 2y e^{-y/2};
$$
while when
$O_\ell -1 \geq e^{y/2}$, which is equivalent to $y\leq 2 \log(O_\ell -1)$, we have
$$
2y \left[1-(1-e^{-y})^{O_\ell -1}\right] \leq 2y \leq 4\log(O_\ell -1).
$$
Hence, combining \eqref{Ineq-7} and \eqref{Ineq-8}, we get
\begin{equation}\label{Ineq-9}
\Q^{Z}\left[ \log_+^2 \overline{W}_{s-\Gamma_\ell}^{Z,\Gamma_\ell} \Big|w, \Gamma_\ell, O_\ell  \right]
\leq 18\log^2(O_\ell -1)
 + \int_0^\infty 4y e^{-y/2}\mathrm{d} y .
\end{equation}
By \eqref{Ineq-10}, \eqref{Ineq-6} and \eqref{Ineq-9}, we obtain
\begin{align}\label{Ineq-11}
\Q^{Z}\left[\log_+^2(3H_2)\right] \leq & 2(\log 3)^2 + 2\lambda_0^2 \psi'(\lambda^*)m s^2 + 4\Q^{Z}\left[\sum_{\ell \geq 1} 1_{\{\Gamma_\ell \leq s\}} 18 \log^2 (O_\ell -1)   \right]\nonumber\\
&+ 4\int_0^\infty 4y e^{-y/2}\mathrm{d} y \Q^{Z}\left[\sum_{\ell \geq 1} 1_{\{\Gamma_\ell \leq s\}} \right]
=
K_1 + K_2 s + K_3 s^2,
\end{align}
here
\begin{align*}
K_1 & = 2(\log 3)^2 ,\quad K_2=  4\psi'(\lambda^*) m \int_0^\infty 4y e^{-y/2}\mathrm{d} y + 72\psi'(\lambda^*)  \sum_{k\geq 2}k \log^2(k-1) p_k ,\\
K_3 & = 2\lambda_0^2 \psi'(\lambda^*) m.
\end{align*}
By \eqref{Equal-0}, \eqref{Ineq-3}, \eqref{Ineq-4}, \eqref{Ineq-5} and \eqref{Ineq-11}, we deduce that $\mathbf{P}_{(\delta_0, \delta_0)}\left[W_1^{\mathcal{Z}}\log_+^2 W_1^{\mathcal{Z}}\right] < \infty$.
		
$Step\ 3:$
In this step we prove the second inequality of \eqref{Check2}. We use similar arguments as in Step 2. First we have
\begin{equation}\label{Equal-2}
\mathbf{P}_{(\delta_0, \delta_0)}\left[D_1^{\mathcal{Z},+}\log_+ D_1^{\mathcal{Z},+}\right] = \int_0^\infty \kappa e^{-\kappa s}\mathrm{d}s \mathbf{P}_{(\delta_0, \delta_0)}\left[D_s^{Z,+}\log_+ D_s^{Z,+}\right],
\end{equation}
here
$$
D_s^{Z,+}:= \lambda_0\langle (\cdot+\lambda_0 s)_+ e^{-\lambda_0(\cdot+\lambda_0 s)}, Z_s \rangle.
$$
For any $\epsilon>0$, there exists a constant $K_\epsilon>0$ such that $\sup_{x\in\R}\left[ (x)_+ e^{-\epsilon x}\right] \leq K_\epsilon.$
Using the definition \eqref{def-add-mart} of the additive martingale $W_t^{Z}(\lambda)$, one can easily get that
$$
D_s^{Z,+} \leq K_\epsilon \lambda_0 \langle e^{-(\lambda_0 -\epsilon)(\cdot + \lambda_0 s)}, Z_s\rangle =  K_\epsilon\lambda_0 e^{\epsilon^2 s/2}W_s^{Z}(\lambda_0 -\epsilon).
$$
By the inequality $\log_+(xy)\leq \log_+ x + \log_+ y$
and the equality $\mathbf{P}_{(\delta_0,\delta_0)}\left[W_s^Z(\lambda_0 -\epsilon)\right]=1$,
we get
\begin{align}\label{Ineq-12}
&\mathbf{P}_{(\delta_0, \delta_0)}\left[D_s^{Z,+}\log_+ D_s^{Z,+}\right] \nonumber\\
&\leq  K_\epsilon \lambda_0 e^{\epsilon^2 s/2}\log_+ \left(K_\epsilon \lambda_0 e^{\epsilon^2 s/2}\right) +K_\epsilon \lambda_0 e^{\epsilon^2 s/2}\mathbf{P}_{(\delta_0, \delta_0)}\left[W_s^{Z}(\lambda_0 -\epsilon)\log_+ W_s^{Z}(\lambda_0 -\epsilon)\right] .
\end{align}
By \eqref{Equal-2} and \eqref{Ineq-12}, to complete the proof, it suffices to prove that, for fixed $\epsilon^2/2 < \kappa $, we have
\begin{equation}\label{Next-Goal}
\int_0^\infty e^{-(\kappa -\epsilon^2 /2)s}\mathrm{d}s \mathbf{P}_{(\delta_0, \delta_0)}\left[W_s^{Z}(\lambda_0 -\epsilon)\log_+ W_s^{Z}(\lambda_0 -\epsilon)\right] < \infty.
\end{equation}
As in Step\ 2, we define $\Q^{Z,\epsilon}$ by
$$
\frac{\mathrm{d} \Q^{Z,\epsilon}}{\mathrm{d}\mathbf{P}_{(\delta_0, \delta_0)}}\bigg|_{\sigma(Z_r, r\leq s)} := W_s^{Z}(\lambda_0 -\epsilon),\quad s\geq 0.
$$
Then $Z$ has another spine decomposition, which is the same as the spine decomposition at the beginning of Step 2 except with $\lambda_0$ replaced by $\lambda_0-\epsilon$, also see \cite[page 59--60]{Ky}.
Set $g(t)= e^{-\epsilon^2 t/2 -(\lambda_0-\epsilon)\lambda_0 t}.$ Using the same notation as in Step 2,  we have
\begin{align*}
W_s^{Z}(\lambda_0 -\epsilon) = &\sum_{\ell\geq 1} 1_{\{\Gamma_\ell \leq s  \}} g(\Gamma_\ell) W_{s-\Gamma_\ell}^{Z,\Gamma_\ell}(\lambda_0 -\epsilon) 1_{\left\{g(\Gamma_\ell) W_{s-\Gamma_\ell}^{Z,\Gamma_\ell}(\lambda_0 -\epsilon)<A\right\}}\\ &  +\sum_{\ell\geq 1} 1_{\{\Gamma_\ell \leq s  \}}g(\Gamma_\ell) W_{s-\Gamma_\ell}^{Z,\Gamma_\ell}(\lambda_0 -\epsilon) 1_{\left\{g(\Gamma_\ell) W_{s-\Gamma_\ell}^{Z,\Gamma_\ell}(\lambda_0 -\epsilon)\geq A\right\}}+ g(s)e^{-(\lambda_0 -\epsilon)w_s}\\ =: &  H_1 + H_2 + H_3 ,
\end{align*}
where $A>1$ is a constant such that
$
\log A  > 1 \geq \sup_{a\geq 1} \left[\log(1+a) - \log a\right],
$
which means that $\log(b+c)\leq \log b + \log c$ for  all $b,c \geq A$.
Also note that \eqref{Ineq-3} and $H_1 \leq A\sum_{\ell\geq 1}1_{\{\Gamma_\ell \leq s\}}$ still hold. And we have
\begin{align*}
\Q^{Z,\epsilon}[\log_+(3H_3) ] &\leq \log 3 + s\epsilon(\lambda_0 -\epsilon/2) + (\lambda_0 -\epsilon)\Q^{Z,\epsilon}\vert w_s + (\lambda_0 -\epsilon)s\vert\\
& = \log 3 + s\epsilon(\lambda_0 -\epsilon/2) + (\lambda_0 -\epsilon)\sqrt{\frac{2}{\pi}}\sqrt{s}.
\end{align*}
Similarly we define $\overline{W}_{s-\Gamma_\ell}^{Z,\Gamma_\ell}(\lambda_0 -\epsilon)$ by
$$
\overline{W}_{s-\Gamma_\ell}^{Z,\Gamma_\ell}(\lambda_0 -\epsilon) := e^{(\lambda_0-\epsilon)w_{\Gamma_\ell} } {W}_{s-\Gamma_\ell}^{Z,\Gamma_\ell}(\lambda_0 -\epsilon).
$$
Then using an argument similar to \eqref{Ineq-10}, we have
\begin{align*}
&\log_+ (3H_2)  \leq \log 3 + \log_+ H_2 \\  &\leq \log 3+ \sum_{\ell \geq 1} 1_{\{\Gamma_\ell \leq s  \}} \log_+\left(g(\Gamma_\ell)e^{-(\lambda_0-\epsilon) w_{\Gamma_\ell}}\right) + \sum_{\ell \geq 1} 1_{\{\Gamma_\ell \leq s  \}} \log_+ \overline{W}_{s-\Gamma_\ell}^{Z,\Gamma_\ell}(\lambda_0 -\epsilon)
\end{align*}
and
$$ \Q^{Z,\epsilon}\left[\sum_{\ell \geq 1} 1_{\{\Gamma_\ell \leq s  \}} \log_+\left(g(\Gamma_\ell)e^{-(\lambda_0-\epsilon) w_{\Gamma_\ell}}\right) \right]\leq \psi'(\lambda^*)m \int_0^s \left[\Q^{Z,\epsilon}\left|w_r+(\lambda_0-\epsilon)r\right|+\epsilon\left(\lambda_0-\frac{\epsilon}{2}\right)r\right] \mathrm{d}r.$$
Since \eqref{Ineq-7} and \eqref{Ineq-8} hold with $\overline{W}_{s-\Gamma_\ell}^{Z,\Gamma_\ell}$ replaced by $\overline{W}_{s-\Gamma_\ell}^{Z,\Gamma_\ell}(\lambda_0 -\epsilon)$ (we only use the martingale property and branching property),
\eqref{Ineq-9} holds for $\overline{W}_{s-\Gamma_\ell}^{Z,\Gamma_\ell}(\lambda_0 -\epsilon)$. Applying Jensen's inequality for $\overline{W}_{s-\Gamma_\ell}^{Z,\Gamma_\ell}(\lambda_0 -\epsilon)$ in \eqref{Ineq-9}, we finally deduce that there exist constants $K_j^\epsilon, j=1,2,3,4,5$, such that
for all $s\geq 0$,
\begin{equation}\label{Ineq-13}
\mathbf{P}_{(\delta_0, \delta_0)}\left[W_s^{Z}(\lambda_0 -\epsilon)\log_+ W_s^{Z}(\lambda_0 -\epsilon)\right] \leq K_1^\epsilon + K_2^\epsilon \sqrt{s} + K_3^\epsilon s + K_4^\epsilon s^{3/2} + K_5^\epsilon s^2.
\end{equation}
Combining \eqref{Ineq-12}, \eqref{Next-Goal} and \eqref{Ineq-13},
we obtain
$\mathbf{P}_{(\delta_0, \delta_0)}\left[D_1^{\mathcal{Z},+}\log_+ D_1^{\mathcal{Z},+}\right]<\infty.$
\hfill$\Box$
		
\begin{lemma}\label{Integrability-condition2}
If \eqref{condition2} holds, then $\sum_{n\geq 1} n(\log n)^2 p_n < \infty.$
\end{lemma}
\textbf{Proof :} By the definition of $\{p_n : n\geq 2\}$, we only need to prove that
\begin{equation}\label{Ineq-14}
\int_{(0,\infty)} \sum_{n\geq2} n(
\log n)^2 \frac{(\lambda^* x)^n}{n!} e^{-\lambda^* x}\nu(\mathrm{d} x) < \infty.
\end{equation}
Define $h(x):= (\log(1+x))^2,$ then
$h''(x) = \frac{2}{(1+x)^2}\left(1-\log(1+x)\right)$. When $x \geq 2> e-1, h''(x) < 0$, which implies $h$ is concave in $[2,\infty).$ By Jensen's inequality,
\begin{align}\label{Ineq-15}
&\sum_{n\geq3} n(\log n)^2 \frac{(\lambda^* x)^n}{n!} e^{-\lambda^* x}  = \lambda^* x \sum_{n\geq2} \left(
\log (1+n)\right)^2 \frac{(\lambda^* x)^n}{n!} e^{-\lambda^* x} \nonumber\\
&\leq  (\lambda^* x) \left[\sum_{n\geq 2}\frac{(\lambda^* x)^n}{n!} e^{-\lambda^* x}  \right]\left\{\log\left[\frac{\sum_{n\geq 2} n(\lambda^* x)^n e^{-\lambda^* x}/n!}{\sum_{n\geq 2} (\lambda^* x)^n e^{-\lambda^* x}/n!} +1\right] \right\}^2\nonumber\\
&\leq  \lambda^* x \left\{\log\left[\frac{\lambda^* x(1-e^{-\lambda^* x})}{1-e^{-\lambda^* x } -e^{-\lambda^* x}\lambda^* x} +1\right] \right\}^2.
\end{align}
Since
$$\lim_{x\to\infty} \log\left[\frac{\lambda^* x(1-e^{-\lambda^* x})}{1-e^{-\lambda^* x } -e^{-\lambda^* x}\lambda^* x} +1\right] / \log x =1,$$ there exists $K> 0$ such that when $x\geq K,$ we have
\begin{equation}\label{Ineq-16}
\log\left[\frac{\lambda^* x(1-e^{-\lambda^* x})}{1-e^{-\lambda^* x } -e^{-\lambda^* x}\lambda^* x} +1\right]\leq 2 \log x.
\end{equation}
Together with \eqref{Ineq-14}, \eqref{Ineq-15} and \eqref{Ineq-16}, we complete the proof.
\hfill$\Box$

\textbf{Proof of Theorem \ref{th2}: }
By the first two  paragraphs of this section, to prove Theorem \ref{th2}, it suffices to show that,
the limsup in \eqref{e:th2} is valid $\mathbf{P}_{(\delta_0, \delta_0)}$-almost surely.
		
$Case\ 1:\beta \neq 0$.
Let $L_t^{Z}$ be the left-most point of $Z_t$.
Suppose that the times of the continuous immigrations
in the skeleton decomposition of $X$ along the trajectory of $L_t^{Z}$ are given by $\left\{(\tau_n, \bar{X}^{(1,\tau_n)}): n=1,2,...\right\}$,
then it is obvious that $\left\{\tau_n-\tau_{n-1}: n =1,2,... \right\}$ are iid and independent of $Z$, also the law of $\tau_n-\tau_{n-1}$ is exponential with parameter $\kappa= 2\beta.$
		
Since \eqref{condition2} holds,
using Lemmas \ref{Integrability-condition1} and \ref{Integrability-condition2} with $T_n=\tau_n$,
we know that $\mathcal{Z}_n$ satisfies \eqref{moment1}, \eqref{moment2} and \eqref{moment3}. Note that the left support of $\mathcal{Z}_n$ is $\lambda_0 (L_{\tau_n}^{Z} +\lambda_0 \tau_n)$, by \cite[Theorem 6.1]{AESZ},
\begin{equation}\label{Liminf_1}
\liminf_{n\to\infty} \left(\lambda_0 (L_{\tau_n}^{Z} +\lambda_0 \tau_n)-\frac{1}{2}\log n \right)= -\infty,\quad \mathbf{P}_{(\delta_0, \delta_0)}\text{-a.s.}
\end{equation}
By the strong law of large numbers, $\tau_n/n \to (2\beta)^{-1}$ as $n\to\infty$. Hence, \eqref{Liminf_1} is equivalent to
\begin{equation}\label{Liminf_2}
\liminf_{n\to\infty} \left(\lambda_0 (L_{\tau_n}^{Z} +\lambda_0 \tau_n)-\frac{1}{2}\log \tau_n \right)= -\infty,\quad \mathbf{P}_{(\delta_0, \delta_0)}\text{-a.s.}
\end{equation}
Define $W_t^{\Lambda}$ by
$$
W_t^{\Lambda}:= \langle e^{-\lambda_0 (\cdot +\lambda_0 t)}, \Lambda_t\rangle,$$
then
\begin{equation}\label{STEP-1}
\sqrt{\tau_{n} +1} \langle e^{-\lambda_0(\cdot + \lambda_0 (\tau_{n}+1))}, \Lambda_{\tau_{n}+1}\rangle \geq \sqrt{\tau_{n}}  \langle e^{-\lambda_0(\cdot + \lambda_0  (\tau_{n}+1))}, \bar{X}^{(1,\tau_n)}_{1}\rangle= : H_n J_n.
\end{equation}
Here $H_n$ and $J_n$ are defined as
$$
H_n:= \sqrt{\tau_n}e^{-\lambda_0(L_{\tau_n}^{Z} +\lambda_0 \tau_n)},\quad
J_n:=e^{-\lambda_0^2} \langle e^{-\lambda_0(\cdot -L_{\tau_n}^{Z})}, \bar{X}^{(1,\tau_n)}_{1}\rangle.
$$
Then by the construction of the continuous immigration in the skeleton decomposition and the spatial homogeneity of super-Brownian motion, we deduce that $\{J_n:n=1,2,...\}$ are iid and for every $n, J_n$ is independent to $\sigma(H_\ell, \ell \geq 1).$
Define $\mathcal{G}_n:= \sigma(H_\ell, J_\ell: 1\leq \ell \leq n).$  By \eqref{Liminf_2}, we have $\limsup_{n\to\infty} H_n = +\infty$, $\mathbf{P}_{(\delta_0, \delta_0)}$-a.s., which together with  the second Borel-Cantelli lemma (see e.g. \cite[Theorem 5.3.2]{Durrett}) is equivalent to that, for any $K>0$,
\begin{equation}\label{Sum-is-Infinite}
\sum_{n=1}^\infty \mathbf{P}_{(\delta_0, \delta_0)}\left[H_n>K \big| \mathcal{G}_{n-1}\right]= +\infty,\quad \mathbf{P}_{(\delta_0, \delta_0)}\text{-a.s.}
\end{equation}
Now it is clear that $\mathbf{P}_{(\delta_0, \delta_0)}(J_n > 0)>0$, so there exists a constant $\varepsilon>0$ such that for all $n\geq 1, \mathbf{P}_{(\delta_0, \delta_0)}(J_n > \varepsilon)>0$. By \eqref{Sum-is-Infinite} and the independence between $J_n$ and $\mathcal{G}_{n-1}$, we deduce that, for any $K>0$,
\begin{align*}
&\sum_{n=1}^\infty \mathbf{P}_{(\delta_0, \delta_0)}\left[H_n J_n>K \big| \mathcal{G}_{n-1}\right] \geq 	\sum_{n=1}^\infty \mathbf{P}_{(\delta_0, \delta_0)}\left[J_n> \varepsilon,\ H_n >K /\varepsilon\big| \mathcal{G}_{n-1}\right]\\ &=\mathbf{P}_{(\delta_0, \delta_0)}[J_1> \varepsilon]\sum_{n=1}^\infty \mathbf{P}_{(\delta_0, \delta_0)}\left[H_n >K /\varepsilon\big| \mathcal{G}_{n-1}\right] =+\infty.\quad \mathbf{P}_{(\delta_0, \delta_0)}\text{-a.s.},
\end{align*}
which is, according to the second Borel-Cantelli lemma, equivalent to
\begin{equation}\label{Limsup_1}
\limsup_{n\to\infty} H_nJ_n = +\infty,\quad \mathbf{P}_{(\delta_0, \delta_0)}\text{-a.s.}
\end{equation}
In view of \eqref{STEP-1} and \eqref{Limsup_1}, we get
$$
\limsup_{t\to\infty} \sqrt{t} W_t^\Lambda \geq \limsup_{n\to\infty}\sqrt{\tau_{n} +1} \langle e^{-\lambda_0(\cdot + \lambda_0 (\tau_{n}+1))}, \Lambda_{\tau_{n}+1}\rangle = +\infty,\quad  \mathbf{P}_{(\delta_0, \delta_0)}\text{-a.s.},
$$
which implies the desired result.
	
$Case\ 2: \nu\neq 0$. Suppose that $\nu\left((\varepsilon, +\infty)\right)>0$, then $\nu\left((\varepsilon, +\infty)\right) < \infty.$
Suppose that the times and masses of the discrete immigration along the trajectory of $L_t^{Z}$ in the skeleton decomposition with initial immigration mass large than $\varepsilon$ are
$\left\{(\tilde{\tau}_n, \mathfrak{m}_n): n=1,2,...\right\}$, then
$\{\tilde{\tau}_n -\tilde{\tau}_{n-1}: n=1,2,...\}$ are iid exponential random variables with parameter $\kappa=\int_{(\varepsilon , \infty)} ye^{-\lambda^* y} \nu(\mathrm{d}y)$,  $\mathfrak{m}_n> \varepsilon$ for all $n\geq 1$ with law $ye^{-\lambda^* y} 1_{\{y >\varepsilon\}}\nu(\mathrm{d}y)/ \int_{(\varepsilon , \infty)} ye^{-\lambda^* y} \nu(\mathrm{d}y)$, and $\{\tilde{\tau}_n: n=1,2,...\}$ is independent of  $Z$.
Applying Lemmas \ref{Integrability-condition1} and \ref{Integrability-condition2}, with $T_n=\tilde{\tau}_n$, we get
\begin{equation}\label{Liminf_3}
\liminf_{n\to\infty} \left(\lambda_0 (L_{\tilde{\tau}_n}^{Z} +\lambda_0 \tilde{\tau}_n)-\frac{1}{2}\log \tilde{\tau}_n \right)= -\infty,\quad \mathbf{P}_{(\delta_0, \delta_0)}\text{-a.s.}
\end{equation}
By the same argument as Case 1, we have
\begin{equation}\label{STEP-2}
\sqrt{\tilde{\tau}_{n}} \langle e^{-\lambda_0(\cdot + \lambda_0 \tilde{\tau}_{n})}, \Lambda_{\tilde{\tau}_{n}}\rangle \geq \sqrt{\tilde{\tau}_{n}} e^{-\lambda_0 (L_{\tilde{\tau}_n}^{Z,1} +\lambda_0 \tilde{\tau}_n)}\mathfrak{m}_n > \varepsilon \sqrt{\tilde{\tau}_{n}} e^{-\lambda_0 (L_{\tilde{\tau}_n}^{Z,1} +\lambda_0 \tilde{\tau}_n)} .
\end{equation}
Combining \eqref{Liminf_3} and \eqref{STEP-2}, we also get the desired results.
\hfill
$\Box$

A byproduct of the proof of Theorem \ref{th2} is the following result:

\begin{cor}\label{Cor1}
Let $L_t$ be the minimum of the support of $X_t$, i.e., $L_t:= \inf\{y\in \R: X_t\left((-\infty, y)\right)> 0\}$. If \eqref{condition1} and \eqref{condition2} hold, then on $\mathcal{E}^c$, it holds that
\begin{equation}\label{Liminf_0}
\liminf_{t\to\infty} \left(L_t +\lambda_0 t -\frac{1}{2\lambda_0} \log t\right)= -\infty \quad \P\mbox{-almost surely}.
\end{equation}
\end{cor}
\textbf{Proof:}
Let $L_t^{\Lambda}$ be the minimum of the support of $\Lambda_t$. We keep
the notation in the proof of Theorem \ref{th2}.

If $\nu \neq 0$,
by the definition of $L_{\tilde{\tau}_{n}}^{\Lambda}$, we  have
$L_{\tilde{\tau}_{n}}^{\Lambda}\leq L_{\tilde{\tau}_n}^{Z}$, $\forall n\ge 1$, $\mathbf{P}_{(\delta_0,\delta_0)}$-a.s.
By the branching property, we deduce that on $\left(\mathcal{E}^{\Lambda}\right)^c$,
$L_{\tilde{\tau}_{n}}^{\Lambda}\leq L_{\tilde{\tau}_n}^{Z}$,  $\forall n\ge 1$, $\mathbf{P}_{(\delta_0,\delta_0)}$-a.s.
Together with \eqref{Liminf_3}, we get \eqref{Liminf_0}.

If  $\beta \neq 0$,
for a fixed constant $A$, define $\mathcal{J}_n$ by $$\mathcal{J}_n:= \langle 1_{(-\infty, A+ L_{\tau_n}^{Z})}(\cdot), \bar{X}^{(1,\tau_n)}_{1}\rangle = \langle 1_{(-\infty, A)}(\cdot-L_{\tau_n}^{Z}), \bar{X}^{(1,\tau_n)}_{1}\rangle.$$
Put $\mathcal{H}_n:= \lambda_0 (L_{\tau_n}^{Z} +\lambda_0 \tau_n)-\frac{1}{2}\log \tau_n.$
By the spatial homogeneity of
super-Brownian motion, $\{\mathcal{J}_n\}$ are iid and for every $n, \mathcal{J}_n$ is independent of $\sigma(\mathcal{H}_\ell , \ell \geq 1)$. We also define $\widetilde{\mathcal{G}}_n := \sigma(\mathcal{H}_\ell, \mathcal{J}_\ell, 1\leq \ell \leq n)$.
Since $\mathbf{P}_{(\delta_0,\delta_0)}\left(\Vert \bar{X}^{(1,\tau_n)}_{1} \Vert > 0\right)= \mathbf{P}_{(\delta_0,\delta_0)}\left(\Vert \bar{X}^{(1,\tau_1)}_{1} \Vert > 0\right)>0$
and $\lim_{A\to+\infty} \mathcal{J}_n = \Vert \bar{X}^{(1,\tau_n)}_{1} \Vert, \mathbf{P}_{(\delta_0,\delta_0)}$-a.s.,
there exists an $A$ such that
$\mathbf{P}_{(\delta_0,\delta_0)}(\mathcal{J}_n >0 )=\mathbf{P}_{(\delta_0,\delta_0)}(\mathcal{J}_1 >0 )>0$.
We see that for any $K>0$,
\begin{align*}
\sum_{n=1}^\infty \mathbf{P}_{(\delta_0, \delta_0)}\left[\mathcal{J}_n > 0, \mathcal{H}_n <-K \big| \widetilde{\mathcal{G}}_{n-1}\right]
=  \mathbf{P}_{(\delta_0, \delta_0)}\left[\mathcal{J}_1 > 0\right]\sum_{n=1}^\infty \mathbf{P}_{(\delta_0, \delta_0)}\left[ \mathcal{H}_n <-K \big| \widetilde{\mathcal{G}}_{n-1}\right] = +\infty,
\end{align*}
$\mathbf{P}_{(\delta_0, \delta_0)}\mbox{-a.s.}$, where in the last equality we used \eqref{Liminf_2} and the second Borel-Cantelli lemma.
Therefore, for all $K>0, \mathbf{P}_{(\delta_0, \delta_0)}\left(\mathcal{J}_n > 0, \mathcal{H}_n <-K\mbox{ i.o.} \right)=1$. Note that $$\left\{\mathcal{J}_n > 0, \mathcal{H}_n <-K \right\}\subset \left\{\lambda_0 (L_{\tau_n+1}^{\Lambda} +\lambda_0 \tau_n)-\frac{1}{2}\log \tau_n<-K + \lambda_0 A\right\},$$
we get $$ \mathbf{P}_{(\delta_0, \delta_0)}\left(\lambda_0 (L_{\tau_n+1}^{\Lambda} +\lambda_0 \tau_n)-\frac{1}{2}\log \tau_n<-K + \lambda_0 A \mbox{ i.o.} \right)=1.$$
Since $(\tau_n +1)/\tau_n \to 1$ as $n\to\infty$ and $K$ is arbitrary, we get that \eqref{Liminf_0} holds $\mathbf{P}_{(\delta_0,\delta_0)}$-almost surely. By the branching property argument, we get the desired result.
\hfill$\Box$

\bigskip
\noindent
{\bf Acknowledgment:}
Part of the research for this paper was done while the third-named author was visiting
Jiangsu Normal University, where he was partially supported by  a grant from
the National Natural Science Foundation of China (11931004) and by
the Priority Academic Program Development of Jiangsu Higher Education Institutions.
\bigskip
\noindent

\begin{singlespace}

\end{singlespace}

\vskip 0.2truein
\vskip 0.2truein

\noindent{\bf Haojie Hou:}  School of Mathematical Sciences, Peking
University,   Beijing, 100871, P.R. China. Email: {\texttt
houhaojie@pku.edu.cn}

\smallskip

\noindent{\bf Yan-Xia Ren:} LMAM School of Mathematical Sciences \& Center for
Statistical Science, Peking
University,  Beijing, 100871, P.R. China. Email: {\texttt
yxren@math.pku.edu.cn}

\smallskip
\noindent {\bf Renming Song:} Department of Mathematics,
University of Illinois,
Urbana, IL 61801, U.S.A.
Email: {\texttt rsong@math.uiuc.edu}

\end{doublespace}


\small
\begin{thebibliography}{99}


\bibitem{AESZ}
A\"{i}d\'ekon, E. and Shi, Z.:
The Seneta-Heyde scaling for the branching random walk.
\emph{Ann. Probab.} \textbf{42} (2014), 959--993.

\bibitem{BKM} Berestycki, J., Kyprianou, A. E. and Murillo-Salas, A.: The prolific backbone for supercritical
superprocesses. \emph{Stoch. Proc. Appl.} \textbf{121} (2011), 1315--1331.

\bibitem{BK1}
Biggins, J. D. and Kyprianou, A. E.:
Branching random walk: Seneta-Heyde norming.
In \emph{Trees (Versailles, 1995) (B. Chauvin et al., eds.),
Progr. Probab.}
\textbf{40} (1996), 31--49, Birkh\"{a}user, Basel.


\bibitem{BK2}
Biggins, J. D. and Kyprianou, A. E.:
Seneta-Heyde norming in the branching
random walk. \emph{Ann. Probab.} \textbf{25} (1997), 337--360.

\bibitem{BK3}
Biggins, J. D. and Kyprianou, A. E.:
Measure change in multitype branching. \emph{Adv. in Appl. Probab.} \textbf{36} (2004), 544--581.


\bibitem{Ch}
Chauvin, B.: Multiplicative martingales and stopping lines for branching Brownian motion,  \emph{Ann. Probab.} \textbf{30} (1991) 1195--1205.


\bibitem{Durrett}
Durrett, R. T.:
\emph{Probability: Theory and Examples.} Fourth edition. Cambridge Series in Statistical and Probabilistic Mathematics, 31. Cambridge University Press, Cambridge, 2010.



\bibitem {E.B1.}
Dynkin, E. B.:
Branching exit Markov systems and superprocesses. \emph{Ann. Probab.} \textbf{29} (2001), 1833--1858.



\bibitem{D}
Dynkin, E. B.:
         Branching particle systems and superprocesses, {\em  Ann. Probab.} {\bf 19(3)} (1991), 1157--1194.

\bibitem{Dyn1993}
Dynkin, E. B.:
Superprocesses and partial differential equations. \emph{ Ann. Probab.} \textbf{21} (1993), 1185--1262.



\bibitem{Dyn2002}
Dynkin, E. B.:
Diffusions, superdiffusions and partial differential equations.
AMS (2002), Providencem R.I.

\bibitem{DyKu} Dynkin, E. B. and Kuznetsov, S. E.: $\mathbb{N}$-measures for branching exit Markov systems and their applications to differential equations.
    \emph{Probab. Theory Relat. Fields.} \textbf{130}(1) (2004), 135--150.

\bibitem{EKW} Eckhoff, M., Kyprianou, A. E. and Winkel, M.: Spine, skeletons and the strong law of large numbers. \emph{Ann.
Probab.}, \textbf{43} (2015), 2594--2659.


\bibitem{HLZ}
He, H., Liu, J.-N. and Zhang, M.:
On Seneta-Heyde scaling for a stable branching random walk. \emph{Adv. in Appl. Probab.} \textbf{50(2)} (2018), 565--599.

\bibitem{Heyde}
Heyde, C. C.:
Extension of a result of Seneta for the super-critical Galton-Watson process. \emph{Ann. Math. Statist.} \textbf{41} (1970), 739--742.



\bibitem{HS}
Hu, Y. and Shi, Z.:
Minimal position and critical martingale convergence in branching random walks, and directed polymers on disordered trees. \emph{Ann. Probab.} \textbf{37} (2009), 742--789.


\bibitem{Im} Imhof, J.-P.:
Density factorizations for Brownian motion, meander and the three dimensional Bessel process, and applications. \emph{J. Appl. Probab.}  \textbf{21} (1984), 500--510.


\bibitem{Ky} Kyprianou, A. E.: Travelling wave solutions to the K-P-P equation: alternatives to Simon Harris' probabilistic analysis.
\emph{ Ann. Inst. H. Poincar\'{e} Probab. Statist.}
\textbf{40} (2004), 53--72.	
   	

\bibitem{KLMR} Kyprianou, A.E., Liu, R.-L., Murillo-Salas, A. and Ren, Y.-X.: Supercritical super-Brownian
motion with a general branching mechanism and travelling waves.
\emph{ Ann. Inst. H. Poincar\'{e} Probab. Statist.}
\textbf{48} (2012), 661--687.


\bibitem{LZ} Li, Z.:
Measure valued branching Markov processes.
Springer, Berlin, 2011.



\bibitem{MP} Maillard, P. and Pain, M.:
1-stable fluctuations in branching Brownian motion at critical temperature I: The derivative martingale.
\emph{Ann. Probab.} \textbf{47}(5) (2019), 2953--3002.


\bibitem{Ne} Neveu, J.:  Multiplicative martingales for spatial branching processes. In \emph{Seminar on Stochastic
Processes}, 1987,
\emph{Progr. Probab.  Statist.},
\textbf{15} (1988) 223--241,
Birkha\"{u}ser, Boston.

\bibitem{RSZ} Ren, Y.-X., Song, R. and Zhang, R.: The extremal process of super-Brownian motion. \emph{Stoch. Proc. Appl.} \textbf{137}(2021), 1--34.




\bibitem{Sen} Seneta E.: On recent theorems concerning the supercritical Galton-Watson
process. \emph{Ann. Math. Statist.} \textbf{39} (1968), 2098--2102.


\bibitem{ShWa} Shiga, T. and Watanabe, S.: Bessel diffusions as a one-parameter family of diffusion processes.
\emph{Z. Wahrsch. Verw. Gebiete}
\textbf{27} (1973), 37--46.


\bibitem{YR}  Yang, T. and Ren, Y.-X.: Limit theorem for derivative martingale at criticality w.r.t. branching Brownian motion. \emph{Statist. Probab. Lett.} {\bf 81} (2011) 195-200.




\end{thebibliography}
\end{document}